\documentclass[12pt,reqno]{amsart}

\usepackage{amssymb}
\usepackage{bbm}
\usepackage{mathrsfs}
\usepackage[all]{xy}
\usepackage{float}
\usepackage{graphicx}
\usepackage{array}
\usepackage{diagbox}

\usepackage{etoolbox}
\patchcmd{\smallmatrix}{\thickspace}{\kern.15em}{}{}

\newtheorem{theo}{Theorem}[section]
\newtheorem{lem}[theo]{Lemma}
\newtheorem{prop}[theo]{Proposition}
\newtheorem{coro}[theo]{Corollary}

\theoremstyle{definition}
\newtheorem{defi}[theo]{Definition}
\newtheorem{ttt}[theo]{}
\newtheorem{conv}[theo]{Convention}

\newtheorem*{resu}{Results}
\newtheorem*{nota}{Conventions and Notation}

\newtheorem*{proop}{Proof of Proposition~\ref{corr}}

\theoremstyle{remark}
\newtheorem{rem}[theo]{Remark}
\newtheorem{rems}[theo]{Remarks}
\newtheorem{ex}[theo]{Example}

\makeatletter
\newcommand*{\Biggg}[1]{{\hbox{$\left#1\vbox to35\p@{}\right.\n@space$}}}
\makeatother

\newcommand{\pmodulo}[1]{\nobreak\mkern8mu(\textup{mod}\,\,#1)}

\newcommand{\br}{ }
\newcommand{\brr}{, }

\newcommand{\C}{\mathop{\text{\rm C}}\nolimits}
\newcommand{\Gal}{\mathop{\text{\rm Gal}}\nolimits}
\newcommand{\Frob}{\mathop{\text{\rm Frob}}\nolimits}
\newcommand{\NS}{\mathop{\text{\rm NS}}\nolimits}
\newcommand{\Pic}{\mathop{\text{\rm Pic}}\nolimits}
\newcommand{\GL}{\mathop{\text{\rm GL}}\nolimits}
\renewcommand{\O}{\mathop{\text{\rm O}}\nolimits}
\newcommand{\SO}{\mathop{\text{\rm SO}}\nolimits}
\newcommand{\End}{\mathop{\text{\rm End}}\nolimits}
\newcommand{\Aut}{\mathop{\text{\rm Aut}}\nolimits}
\newcommand{\M}{\mathop{\text{\rm M}}\nolimits}
\newcommand{\Tr}{\mathop{\text{\rm Tr}}\nolimits}
\newcommand{\im}{\mathop{\text{\rm im}}\nolimits}
\newcommand{\rk}{\mathop{\text{\rm rk}}\nolimits}
\newcommand{\sgn}{\mathop{\text{\rm sgn}}\nolimits}
\newcommand{\Spec}{\mathop{\text{\rm Spec}}\nolimits}
\newcommand{\Sym}{\mathop{\text{\rm Sym}}\nolimits}
\newcommand{\Kum}{\mathop{\text{\rm Kum}}\nolimits}

\newcommand{\maxi}{{\!\text{\rm max}}}
\newcommand{\Hg}{\text{\rm Hg}}
\newcommand{\J}{\text{\rm J}}
\newcommand{\E}{\text{\rm E}}
\newcommand{\K}{\text{\rm K}}
\newcommand{\N}{\text{\rm N}}
\newcommand{\alg}{\text{\rm alg}}
\newcommand{\id}{\text{\rm id}}
\newcommand{\pr}{\text{\rm pr}}
\newcommand{\cl}{\text{\rm cl}}
\newcommand{\an}{\text{\rm an}}
\newcommand{\et}{\text{\rm \'et}}

\newcommand{\bbC}{{\mathbbm C}}
\newcommand{\bbF}{{\mathbbm F}}
\newcommand{\bbN}{{\mathbbm N}}
\newcommand{\bbQ}{{\mathbbm Q}}
\newcommand{\bbR}{{\mathbbm R}}
\newcommand{\bbZ}{{\mathbbm Z}}

\newcommand{\calC}{{\mathscr{C}}}
\newcommand{\calO}{{\mathscr{O}}}
\newcommand{\frakl}{\mathfrak{l}}
\newcommand{\frakL}{\mathfrak{L}}
\newcommand{\frakM}{\mathfrak{M}}

\newcommand{\frakp}{\mathfrak{p}}
\newcommand{\frakP}{\mathfrak{P}}

\newcommand{\Pb}{{\text{\bf P}}}

\newcommand{\Exterior}{\mathchoice{{\textstyle\bigwedge}}%
 {{\bigwedge}}%
 {{\textstyle\wedge}}%
 {{\scriptstyle\wedge}}}

\renewcommand{\atop}[2]{\genfrac{}{}{0pt}{}{#1}{#2}}

\newcounter{abc}
\newenvironment{abc}{\begin{list}{\rm \alph{abc}) }%
{\usecounter{abc} \leftmargin=0.0pt \labelsep=0.0pt %
\listparindent=0.0pt \labelwidth=0.0pt \parsep=\smallskipamount%
 \itemsep=0.0pt \topsep=0.0pt \partopsep=\smallskipamount}}{\end{list}}

\newcounter{iii}
\newenvironment{iii}{\begin{list}{\rm \roman{iii}) }%
{\usecounter{iii} \leftmargin=0.0pt \labelsep=0.0pt %
\listparindent=0.0pt \labelwidth=0.0pt \parsep=\smallskipamount%
 \itemsep=0.0pt \topsep=0.0pt \partopsep=\smallskipamount}}{\end{list}}

\def\rightend#1#2{{%
 \leavevmode\nobreak\hskip .5em plus 1fil
 \penalty600 \hskip 0pt plus -1filll
 \vadjust{}\nobreak\hskip 0pt plus 1filll%
 #1\parfillskip=#2\relax \par}}

\def\eop{\ifmmode\rule[-22pt]{0pt}{1pt}\ifinner\tag*{$\square$}\else\eqno{\square}\fi\else\rightend{$\square$}{0pt}\fi
}

\thanks{}

\title[On the component group of the algebraic monodromy group]{\boldmath On the component group of the algebraic monodromy group of a
$K3$~surface}

\begin{document}

\author{Andreas-Stephan Elsenhans}

\address{Institut f\"ur Mathematik\\ Universit\"at W\"urzburg\\ Emil-Fischer-Stra\ss e 30\\ D-97074 W\"urzburg\\ Germany}
\email{stephan.elsenhans@mathematik.uni-wuerzburg.de}
\urladdr{https://www.mathematik.uni-wuerzburg.de/institut/personal/elsenhans.html}

\author[J\"org Jahnel]{J\"org Jahnel}

\address{\mbox{Department Mathematik\\ \!Univ.\ \!Siegen\\ \!Walter-Flex-Str.\ \!3\\ \!D-57068 \!Siegen\\ \!Germany}}
\email{jahnel@mathematik.uni-siegen.de}
\urladdr{https://wwwuser.gwdg.de/~jjahnel}

\thanks{}

\date{January~10, 2024}



\begin{abstract}
We provide a lower bound for the number of components of the algebraic monodromy group in the situation of a
$K3$~surface
over a number
field~$k$.
In~the CM case, our bound is sharp. As an application, we describe, in the case of CM, the jump character \cite[Definition~2.4.6]{CEJ} entirely in terms of the endomorphism field and the geometric Picard rank.
\end{abstract}

\maketitle
\thispagestyle{empty}

\section{Introduction}

\begin{ttt}
Let
$X$
be a
$K3$~surface
over a number
field~$k$.
Associated
with~$X$,
for every prime number
$l$,
one has a continuous representation~(\cite[Expos\'e~VI, Proposition~1.2.5]{SGA5} and \cite[Expos\'e~VIII, Th\'eor\`eme~5.2]{SGA4})
$$\varrho_{X,\overline{l}}\colon \Gal(\overline{k}/k) \longrightarrow \GL(T_{\overline{l}}) \,,$$
for
$\smash{T_{\overline{l}} \subset H^2_\et(X_{\overline{k}}, \overline\bbQ_l(1))}$
the transcendental part of the cohomology. The 
representation~$\varrho_{X,\overline{l}}$
gives rise to the {\em algebraic monodromy group\/}
$G_{X,\overline{l}}$
of~$X$,
which is a linear algebraic group
over~$\smash{\overline\bbQ_l}$,
usually disconnected. It~is defined to be the Zariski closure
$\smash{G_{X,\overline{l}} := \overline{\im(\varrho_{X,\overline{l}})}}$
of the image
of~$\varrho_{X,\overline{l}}$
(cf.\ Definition \ref{def_alg_mon}).

For a generic
$K3$~surface,
one has
$\smash{G_{X,\overline{l}} = \O(T_{\overline{l}})}$.
There are, however, cases, in which the jump character
$\tau_X\colon \Gal(\overline{k}/k) \to \{\pm1\}$
(cf.\ Paragraph~\ref{jump_char}) is trivial \cite[Examples~5.5 and~5.6]{EJ22b}, which yields that
$\smash{G_{X,\overline{l}} = \SO(T_{\overline{l}})}$.
It~may happen, too, that the algebraic monodromy group is of positive codimension
in~$\smash{\O(T_{\overline{l}})}$,
but only when
$X$
has a nontrivial endomorphism field
$E \supsetneqq \bbQ$.
I.e., when there is real (RM) or complex multiplication (CM). Then,~for the neutral component of the algebraic monodromy group, one~has
\begin{equation}
\label{Tankeev}
G_{X,\overline{l}}^0 \cong (\C_E(\O(T_{\overline{l}})))^0 \,,
\end{equation}
due to the work of S.\,G.\ Tankeev~\cite{Tan90,Tan95}, together with Yu.\,G.\ Zarhin~\cite[Theorem 2.2.1]{Za}. We~discussed upper bounds for the component group
$\smash{G_{X,\overline{l}}/G_{X,\overline{l}}^0}$
in \cite[Lemma~4.10 and Theorem~4.12]{EJ22b}.
\end{ttt}

\begin{rem}
There~can, of course, be no hope for a general lower bound. Indeed, the representation
$\varrho_{X,\overline{l}}$
induces a homomorphism
$\smash{\underline\varrho_{X,\overline{l}}\colon \Gal(\overline{k}/k) \hookrightarrow G_{X,\overline{l}}/G_{X,\overline{l}}^0}$,
which is surjective. Hence, the component group
$\smash{G_{X,\overline{l}}/G_{X,\overline{l}}^0}$
has a {\em splitting field\/}
$k' \supseteq k$
satisfying the conditions that
$\smash{\varrho_{X,\overline{l}}(\Gal(\overline{k}/k')) \subseteq G_{X,\overline{l}}^0}$
and that the induced homomorphism
$\smash{\Gal(k'/k) \to G_{X,\overline{l}}/G_{X,\overline{l}}^0}$
is bijective. Then~the algebraic monodromy group of the base
extension~$\smash{X_{k'}}$
is
$\smash{G_{X_{k'},\overline{l}} \cong G_{X,\overline{l}}^0}$
and, in particular,~connected.
\end{rem}

\begin{resu}
There is, however, a lower bound for
$K3$~surfaces
over number fields that are linearly disjoint to the endomorphism field
$E$.\medskip

\noindent
{\bf Theorem \ref{main}.}
{\em
Let\/
$X$
be a\/
$K3$~surface
over a number
field\/~$k$
with endomorphism
field\/~$E$
and
let\/~$l$
be any prime~number.

\begin{abc}
\item
If\/~$kE \supseteq k$
is a normal extension then one has a natural surjective homomorphism
\begin{equation}
\label{hom_main}
G_{X,\overline{l}}/G_{X,\overline{l}}^0 \twoheadrightarrow \Gal(kE/k) \,.
\end{equation}
\item
In~any case,
$\smash{\#(G_{X,\overline{l}}/G_{X,\overline{l}}^0)}$
is always a multiple
of\/~$[kE:k]$.%
\end{abc}%
}\medskip

\noindent
The~proof shows that the splitting field
of~$\smash{G_{X,\overline{l}}/G_{X,\overline{l}}^0}$
contains~$kE$.
As~the former is normal
over~$k$,
one might want to conclude that actually there is a natural surjective homomorphism
$\smash{G_{X,\overline{l}}/G_{X,\overline{l}}^0 \twoheadrightarrow \Gal((kE)^{(n)}/k)}$,
for
$(kE)^{(n)}$
the normal closure
of~$kE$
over~$k$.
However,~such an improvement is vacuous, at least~conjecturally.\medskip

\noindent
{\bf Theorem~\ref{normality}.}
{\em
Let\/~$X$
be a\/
$K3$~surface
over a
subfield\/
$K$
of\/~$\bbC$
with endomorphism
field\/~$E$.
Suppose~that the Hodge conjecture holds
for\/~$(X \times X)(\bbC)$.
Then~the com\-pos\-ite field\/
$KE$
is normal
over\/~$K$.
}\medskip

\noindent
The~assumption concerning the Hodge conjecture is fulfilled, for instance, for double covers
of~$\Pb^2$
branched over six lines~\cite[Theorem 2]{Sch}.

It~is fulfilled, as well, in the CM case \cite[Theorem 5.4]{RM}, a case in which the methods developed in order to establish Theorem~\ref{main} are particularly~strong.\medskip

\noindent
{\bf Theorem~\ref{jump_CM_alg}.}
{\em
Let\/~$X$
be a\/
$K3$~surface
over a number
field\/~$k$.
Assume~that the endomorphism
field\/~$E$
of\/~$X$
is a CM~field.

\begin{abc}
\item
Then~there is a natural isomorphism
$$G_{X,\overline{l}}/G_{X,\overline{l}}^0 \stackrel{\cong}{\longrightarrow} \Gal(kE/k) \,.$$
I.e.,~the splitting field
of\/~$\smash{G_{X,\overline{l}}/G_{X,\overline{l}}^0}$
is
exactly\/~$kE$.
\item
The~jump character\/
$\tau_X\colon \Gal(\overline{k}/k) \to \{\pm1\}$
is
$$
\sigma \mapsto \left\{
\begin{array}{rrr}
1 \,,                                                                             & \hspace{.7cm}\textstyle\text{~if\/~} \smash{\frac{22-\rk\Pic X_{\overline{k}}}{[E:\bbQ]}} \text{~is~even} \,, \\
\text{sign~of~the~natural~permutation~}                                           & \\[-.4mm]
\text{action~of~} \sigma \in \smash{\Gal(\overline{k}/k)} \text{~on~the~}         & \\[-.4mm]
\text{conjugates\/~} u_1, \overline{u}_1, \ldots, u_d, \overline{u}_d \text{~of~} & \\[-.4mm]
\text{a~primitive~element~} u_0 \in E \,,                                         & \textstyle\text{~if\/~} \smash{\frac{22-\rk\Pic X_{\overline{k}}}{[E:\bbQ]}} \text{~is~odd} \,. \;\;
\end{array}
\right.
$$
\end{abc}
}

\noindent
Finally, and again conditional to the Hodge conjecture, the real or complex multiplication on
a~$K3$~surface
over any
field~$K \subseteq \bbC$
by some
element~$u \in E$
may be described by a correspondence, which is a
\mbox{$\bbQ$-cycle}
of
codimension~$2$
in~$(X \!\times\! X)_{\overline{K}}$.
As~a digression, we show in Theorem~\ref{corr_bound} that the field of definition of such a correspondence must contain
$K(u)$.
\end{resu}

\begin{rem}
It is our understanding throughout the article that the algebraic monodomy group is a linear algebraic group over the algebraically closed field
$\overline\bbQ_l$.
This is a convention just for convenience. One might as well rely on
\mbox{$l$-adic}
cohomology instead of the
\mbox{$\overline{l}$-adic}
theory. Then one ends up with a group scheme
$G_{X,l}$
defined over
$\bbQ_l$
that underlies the algebraic monodomy group
$\smash{G_{X,\overline{l}}}$.

The component group is, however, exactly the~same, i.e.\
$\smash{G_{X,l}/G_{X,l}^0 \cong G_{X,\overline{l}}/G_{X,\overline{l}}^0}$.
Indeed, by construction, in the irreducible components of
$G_{X,l}$,
the
$\bbQ_l$-rational
points are Zariski dense. Hence, every irreducible component
of~$G_{X,l}$
is in fact geometrically irreducible \cite[Corollaire~4.5.19.3$=$Err${}_{{\rm IV}}$20]{EGAIV}.
\end{rem}

\begin{nota}
We~follow standard conventions and use notation that is standard in Algebra and Algebraic~Geometry. In~particular, and perhaps sometimes slightly deviating from~this,

\begin{iii}
\item
we often work over a base field, which is usually denoted
by~$k$
or~$K$.

When~$K \supseteq k$
is a field extension, we write
$\smash{K^{(n)}}$
for the {\em normal closure\/}
of
$K$
over~$k$.
I.e.,~for the extension field
of~$k$
generated by the
$k$-conjugates
of~$K$.

For an arbitrary field
$K$, we denote
by~$\smash{\overline{K}}$
the algebraic~closure.
\item
For~$B = (v_1, \ldots, v_m)$
a basis of a
\mbox{$K$-vector}
space~$T$,
we denote by
$\M_B^B(f)$
the matrix representing the
\mbox{$K$-linear}
map~$f\colon T \to T$.
Similarly,~for the matrix that represents a
\mbox{$K$-bilinear}
form~$\langle.\,,.\rangle\colon T \times T \to K$,
we write
$\M_B(\langle.\,,.\rangle) := (\langle v_i,v_j\rangle)_{1\leq i,j \leq m}$.
\item
We denote the identity matrix of
size~$m$
by~$E_m$.
\item
\label{splitF}
When a surjective homomorphism
$\varrho\colon \Gal(k^{(n)}/k) \twoheadrightarrow G$
is given onto a finite
group~$G$,
then, by the {\em splitting field\/}
of~$\varrho$
(or~$G$),
we mean the intermediate field of
$k^{(n)}/k$
corresponding to
$\ker\varrho$
under the Galois correspondence. The splitting field
of~$G$
is a Galois extension
of~$k$.
\item
For a linear algebraic group
$G$
over an algebraically closed field, we write
$G^0$
to denote its neutral component.
\item
For a finite dimensional vector space
$T$
over an algebraically closed
field~$K$,
we denote
by~$\GL(T)$
the linear algebraic group, whose
\mbox{$K$-rational}
points are the automorphisms
$f\colon T \to T$.

If~$T$
is equipped with a non-degenerate symmetric bilinear
form~$\langle.,.\rangle$,
we
let~$\O(T)$
be the linear algebraic group, whose
\mbox{$K$-rational}
points are the orthogonal maps
$f\colon T \to T$.
I.e.~those such that
$\langle f(x), f(y) \rangle = \langle x, y \rangle$,
for all
$x,y \in T$.

As~usual, we~put
$\SO(T) := \O(T)^0$.
\item
We~say that a
\mbox{$k$-algebra}~$E$
{\em acts
\mbox{$K$-linearly\/}}
on a
\mbox{$K$-vector}
space~$T$
or that there is a
{\em \mbox{$K$-linear}
action\/}
of~$E$
on~$T$
if
$K \supseteq k$
and a homomorphism
$$E \hookrightarrow \End_K(T)$$
of
\mbox{$k$-algebras}
is provided that respects~units.

In~this situation, we let
$\C_E(\O(T))$
be the {\em centraliser\/}
of~$E$
in~$\O(T)$.
I.e.\ the linear algebraic group, whose
\mbox{$K$-rational}
points are the orthogonal maps
$f\colon T \to T$
that commute with the action
of~$E$.
\item
By a
\mbox{\em $\bbQ$-cycle\/}
on a
scheme~$X$, we mean a formal
\mbox{$\bbQ$-linear}
combination of closed integral subschemes
of~$X$,
each of which is of the same codimension.
\item
For a Galois extension
$k' \supseteq k$
and
$\sigma\in\Gal(k'/k)$,
we let
$\sigma$
act
on~$\Spec k'$
by the morphism
$\sigma\colon \Spec k' \to \Spec k'$
induced via contravariant functoriality from
$\sigma^{-1}\colon k' \to k'$.
This~defines a left
\mbox{$\Gal(k'/k)$-action}
on~$\Spec k$.

Let~$X$
be a
\mbox{$k$-scheme}.
Slightly abusing notation, we again write
$\sigma$
for the auto\-morphism
$\smash{X_{k'} = X \times_{\Spec k} \Spec k' \stackrel{(\id,\sigma)}{\longrightarrow} X \times_{\Spec k} \Spec k' = X_{k'}}$
of~$X_{k'}$
obtained by base~change.
\end{iii}
\end{nota}

\section{Technical prerequisites}

\subsubsection*{Independence
of\/~$l$}\leavevmode\medskip

\noindent
Let us recall a few concepts from J.-P.\ Serre's famous McGill notes~\cite[Ch.~I, \S2]{Se68}.

\begin{defi}
Let
$k$
be a number field and
$\varrho\colon \Gal(\overline{k}/k) \to \GL_n(\overline\bbQ_l)$
a continuous representation, for a certain prime
number~$l$.

\begin{abc}
\item
Then
$\varrho$
is called {\em unramified\/} at a prime
$\frakp$
of~$k$
if
$\varrho(I_\frakP) = \{1\}$,
for
$\frakP$
any extension
of~$\frakp$
to~$\smash{\overline{k}}$
and
$I_\frakP$
the inertia~group.
\item
The representation
$\varrho$
is called {\em rational\/} if there is a finite set
$S$
of primes
of~$k$
such that
\begin{iii}
\item
$\varrho$
is unramified at every prime
$\frakp \not\in S$
and
\item
if
$\frakp \not\in S$
then the characteristic polynomial
$\chi(\varrho(\Frob_\frakp)) \in \overline\bbQ_l[T]$
actually has coefficients
in~$\bbQ$.
\end{iii}
\end{abc}
\end{defi}

\begin{defi}
Let
$k$
be a number field and let
$\smash{\varrho\colon \Gal(\overline{k}/k) \to \GL_n(\overline\bbQ_l)}$
and
$\smash{\varrho'\colon \Gal(\overline{k}/k) \to \GL_n(\overline\bbQ_{l'})}$
be two rational continuous representations, for prime
numbers
$l$
and~$l'$.
Then
$\varrho$
and~$\varrho'$
are said to be {\em compatible\/} if there is a finite set
$S$
of primes
of~$k$
such that

\begin{iii}
\item
$\varrho$
and~$\varrho'$
are both unramified at every prime
$\frakp \not\in S$
and
\item
the characteristic polynomials
$\chi(\varrho(\Frob_\frakp)), \chi(\varrho'(\Frob_\frakp)) \in \bbQ[T]$
coincide for every
$\frakp \not\in S$.
\end{iii}
\end{defi}

\begin{ex}
\label{coh_rep}
Let
$X$
be a smooth projective scheme over a number
field~$k$.
Choose~integers
$i$
and
$j$,
as well as a prime
number~$l$.
Then every element
$\sigma\in\Gal(\overline{k}/k)$
induces an automorphism of
$X_{\overline{k}}$
and hence, by functoriality, an automorphism
$\smash{\varrho^{i,j}_{X,\overline{l}}(\sigma) \in \GL(H^i_\et(X_{\overline{k}}, \overline\bbQ_l(j)))}$.
It~is clear that
$$\varrho^{i,j}_{X,\overline{l}}\colon \Gal(\overline{k}/k) \longrightarrow \GL(H^i_\et(X_{\overline{k}}, \overline\bbQ_l(j)))$$
is a representation.

\begin{iii}
\item
One has that
$\smash{\varrho^{i,j}_{X,\overline{l}}}$
is continuous by \cite[Expos\'e~VI, Proposition~1.2.5]{SGA5}, together with \cite[Expos\'e~VIII, Th\'eor\`eme~5.2]{SGA4}.
\item
Moreover,
$\smash{\varrho^{i,j}_{X,\overline{l}}}$
is unramified at every prime
of~$k$
of residue characteristic
$\smash{\neq \!l}$,
at which
$X$
has good reduction. This~is a consequence of the smooth specialisation theorem~\cite[Expos\'e~XVI, Corollaire~2.2]{SGA4}.
\item
Furthermore,
$\smash{\varrho^{i,j}_{X,\overline{l}}}$
is rational, due to the work of P.\ Deligne on the Weil conjectures \cite[Th\'eor\`eme~(1.6)]{De74}.
\item
Finally,~the representations
$\smash{\varrho^{i,j}_{X,\overline{l}}}$,
for~$X$,
$i$,
and~$j$
fixed, but
$l$
running through the set of all primes numbers, are mutually~compatible.

Indeed,~let
$\frakp$
be a prime
of~$k$
of residue characteristic neither
$l$
nor~$l'$,
at which
$X$
has good reduction. Then, by the smooth specialisation theorem,
$\smash{\chi(\varrho^{i,j}_{X,\overline{l}}(\Frob_\frakp))}$
is the same as the characteristic polynomial
of~$\Frob$
on
$\smash{H^i_\et(X_{\overline\bbF_{\!\frakp}}, \overline\bbQ_l(j))}$.
On the other hand, one has that
$\smash{\chi(\varrho^{i,j}_{X,\overline{l}'}(\Frob_\frakp))}$
coincides with the characteristic polynomial
of~$\Frob$
on
$\smash{H^i_\et(X_{\overline\bbF_{\!\frakp}}, \overline\bbQ_{l'}(j))}$.
\;But~the latter two polynomials agree, as a consequence of the Lefschetz trace formula and the Weil conjectures~\cite[Th\'eor\`eme~(1.6)]{De74}.
\item
If~$X$
is connected and
$\smash{i = 2j = \dim X}$
then
$\smash{H^i_\et(X_{\overline{k}}, \overline\bbQ_l(j))}$
carries a non-degenerate symmetric bilinear form, given by the cup product and Poincar\'e duality \cite[Expos\'e~XVIII, Th\'eor\`eme 3.2.5]{SGA4}. One~then has
$\smash{\im \varrho^{i,j}_{X,\overline{l}} \subseteq \O(H^i_\et(X_{\overline{k}}, \overline\bbQ_l(j)))}$.
\end{iii}
\end{ex}

\begin{ex}
\label{comp}
Let~$X$
be a
$K3$~surface
over a number
field~$k$.
Choose a prime
number~$l$.
For~$c_1$
the first Chern class homomorphism,~put
$$\smash{H_{\overline{l},\alg} := \im(c_1 \!\otimes_\bbZ\! \overline\bbQ_l\colon \Pic X_{\overline{k}} \otimes_\bbZ\! \overline\bbQ_l \to H^2_\et(X_{\overline{k}}, \overline\bbQ_l(1)))}$$
and
$T_{\overline{l}} := (H_{\overline{l},\alg})^\perp$.
Then,~let
$$\varrho_{X,\overline{l}}\colon \Gal(\overline{k}/k) \longrightarrow \O(T_{\overline{l}}) \subset \GL(T_{\overline{l}})$$
be the restriction of
$\smash{\varrho^{2,1}_{X,\overline{l}}}$
to~$\smash{T_{\overline{l}}}$.

\begin{iii}
\item
As~a subrepresentation
of~$\smash{\varrho^{2,1}_{X,\overline{l}}}$,
one immediately has that
$\smash{\varrho_{X,\overline{l}}}$
is continuous. Moreover,
$\smash{\varrho_{X,\overline{l}}}$
is unramified at every prime, at which
$\smash{\varrho^{2,1}_{X,\overline{l}}}$~is.
\item
Furthermore,
$\smash{\varrho_{X,\overline{l}}}$
is rational. The~representations
$\smash{\varrho_{X,\overline{l}}}$,
for
$X$~fixed,
but
$l$
running through the set of all primes numbers, are mutually~compatible.

Indeed,~the characteristic polynomial of the action
on~$\smash{\Pic X_{\overline{k}}}$
of any element
$\smash{\sigma\in\Gal(\overline{k}/k)}$
has rational, in fact integral, coefficients, which trivially do not depend
on~$l$.
The~claim directly follows from~this.
\end{iii}
\end{ex}

\begin{defi}
\label{def_alg_mon}
Let
$k$
be a number field and let\vspace{.1mm}
$\smash{\varrho\colon \Gal(\overline{k}/k) \to \GL_n(\overline\bbQ_l)}$
be a con\-tin\-u\-ous representation. Then~the Zariski closure
$\smash{G_\varrho := \overline{\im(\varrho\,}\!)}$
of the image
of~$\varrho$
is called the {\em algebraic monodromy group\/} of the 
representation~$\varrho$.
It~is clear that
$\smash{G_\varrho}$
is a linear algebraic group
over~$\smash{\overline\bbQ_l}$.
\end{defi}

\begin{rems}
\begin{iii}
\item
Let
$l$
and~$l'$
be two prime numbers and let
$\smash{\varrho\colon \Gal(\overline{k}/k) \to \GL_n(\overline\bbQ_l)}$
and
$\smash{\varrho'\colon \Gal(\overline{k}/k) \to \GL_n(\overline\bbQ_{l'})}$
be compatible rational representations. Then~the algebraic monodromy groups
$G_\varrho$
and~$G_{\varrho'}$
certainly do not coincide, if~only because they are defined over different base~fields.
Concerning~the component groups, there is, however, Theorem~\ref{Serre_Lettre} below, which we are going to make use~of.
\item
The representations considered in Example~\ref{coh_rep} give rise to algebraic monodromy groups. Similarly, for a
$K3$~surface
$X$
over a number
field~$k$
and a prime
number~$l$,
there is the algebraic monodromy group
$\smash{G_{X,\overline{l}} := G_{\varrho_{X,\overline{l}}}}$,
which is called the {\em algebraic monodromy group
of\/}~$X$.
\end{iii}
\end{rems}

\begin{theo}[J.-P.\ Serre]
\label{Serre_Lettre}
Let\/~$k$
be a number field and let\/
$l$
and\/~$l'$
be two prime numbers. Suppose that the two rational representations\/
$\smash{\varrho\colon \Gal(\overline{k}/k) \to \GL_n(\overline\bbQ_l)}$
and\/
$\smash{\varrho'\colon \Gal(\overline{k}/k) \to \GL_n(\overline\bbQ_{l'})}$
are~compatible. Then~the induced homomorphisms
$$\underline\varrho\colon \Gal(\overline{k}/k) \longrightarrow G_\varrho/G_\varrho^0 \qquad\text{and}\qquad \underline\varrho'\colon \Gal(\overline{k}/k) \longrightarrow G_{\varrho'}/G_{\varrho'}^0$$
to the component groups have the same~kernel.\medskip

\noindent
{\bf Proof.}
{\em
This is shown in \cite[Lettre du 29/1/1981, p.~18, Th\'eor\`eme]{Se81}.
}
\eop
\end{theo}

\begin{coro}
\label{split_l}
Let\/~$\varrho_l\colon \Gal(\overline{k}/k) \to \GL_n(\overline\bbQ_l)$,
for\/~$l$
running through the set of all primes numbers, be mutually compatible rational representations. Then the splitting field of the induced homomorphism
$$\underline\varrho_l\colon \Gal(\overline{k}/k) \longrightarrow G_{\varrho_l}/G_{\varrho_l}^0$$
to the component group does not depend
on\/~$l$.
\eop
\end{coro}

\subsubsection*{$l$-adic
cohomology versus algebraic de Rham cohomology}

\begin{ttt}[The comparison isomorphism]
Let~$k_\frakl$
be a local field of
characteristic~$0$
and residue
characteristic~$l>0$.
Then~$\bbQ_l \subseteq k_\frakl$.
We~fix inclusions
$\smash{k_\frakl \subseteq \overline{k}_\frakl = \overline\bbQ_l \subset \bbC_l}$,
for~$\bbC_l$
the completion
of~$\smash{\overline\bbQ_l}$,
which is again algebraically closed~\cite[Theorem~13]{Ko}. Moreover, let
$X$
be a
$K3$~surface
over~$k_\frakl$.
Then, as a particular case of
\mbox{$p$-adic}
Hodge theory \cite[Theorem~III.4.1]{Fa}, there is a natural
$\smash{\Gal(\overline{k}_\frakl/k_\frakl)}$-equivariant
\mbox{$\bbC_l$-linear}
isomorphism
$$
\iota_X\colon H^2_\et(X_{{\overline{k}_\frakl}}, \bbC_l(1)) \stackrel{\cong}{\longrightarrow} \bigoplus_{i=0}^2 H^i(X_{\overline{k}_\frakl},\Omega^{2-i}_{X_{\overline{k}_\frakl}/\overline{k}_\frakl}) \!\otimes_{\overline{k}_\frakl}\! \bbC_l(i\!-\!1)
$$
connecting
\mbox{$l$-adic}
\'etale cohomology with algebraic de Rham cohomology.
\end{ttt}

\begin{ttt}[The transcendental parts]
The comparison
isomorphism~$\iota_X$
is compatible with the first Chern class homomorphisms. I.e., the diagram
$$
\xymatrix{
\Pic X_{\overline{k}_\frakl} \ar@{=}[d] \ar@{->}[rr]^{c_1} && H^2_\et(X_{{\overline{k}_\frakl}}, \bbC_l(1)) \ar@{->}[d]_\cong^{\iota_X} \\
\Pic X_{\overline{k}_\frakl} \ar@{->}[r]^{\!\!\!\!\!\!\!\!\!\!\!\!\!\!\!\!\!\!\!\!\!\!\!\!\!c_1} & H^1(X_{\overline{k}_\frakl},\Omega^1_{X_{\overline{k}_\frakl}/\overline{k}_\frakl}) \,\ar@{^{(}->}[r] & \,\smash{\bigoplus\limits_{i=0}^2} \,H^i(X_{\overline{k}_\frakl},\Omega^{2-i}_{X_{\overline{k}_\frakl}/\overline{k}_\frakl}) \!\otimes_{\overline{k}_\frakl}\! \bbC_l(i\!-\!1)
}
$$
is commutative. Cf.\ \cite[Paragraph~III.4.(a)]{Fa}.
Consequently,~$\iota_X$
induces an~isomorphism
$$
\smash{\iota_X^T\colon T_{\overline{l}} \!\otimes_{\overline\bbQ_l}\! \bbC_l \cong H^0(X_{\overline{k}_\frakl},\Omega^2_{X_{\overline{k}_\frakl}/k_\frakl}) \!\otimes_{k_\frakl}\! \bbC_l(-1) \oplus V_{\overline{k}_\frakl} \!\otimes_{\overline{k}_\frakl}\! \bbC_l \oplus H^2(X_{\overline{k}_\frakl},\calO_{X_{\overline{k}_\frakl}}) \!\otimes_{\overline{k}_\frakl}\! \bbC_l(1) \,,}
$$
for~$\smash{V_{\overline{k}_\frakl} := (c_1(\Pic X_{\overline{k}_\frakl}) \!\otimes_\bbZ\! \overline{k}_\frakl)^\perp \subset H^1(X_{\overline{k}_\frakl},\Omega^1_{X_{\overline{k}_\frakl}/\overline{k}_\frakl})}$.
\end{ttt}

\begin{ttt}[The number field case]
\label{nf_case}
Let~$k$
be a number field and
$X$
a
$K3$~surface
over~$k$.
For~$l$
any prime number, let
$\frakl$
be a prime
of~$k$
lying
above~$l$.

\begin{abc}
\item
Then,~combined with the natural isomorphism
$\smash{H^2_\et(X_{\overline{k}}, \bbC_l(1)) \stackrel{\cong}{\to} H^2_\et(X_{{\overline{k}_\frakl}}, \bbC_l(1))}$
\cite[Expos\'e~XVI, Corollaire~1.6]{SGA4},
$\iota_{X_{k_\frakl}}$
induces a natural~isomorphism
\begin{align*}
\iota_{X,\frakl}\colon H^2_\et(X_{\overline{k}}, \bbC_l(1)) \stackrel{\cong}{\longrightarrow}& \,\bigoplus_{i=0}^2 H^i(X_{\overline{k}_\frakl},\Omega^{2-i}_{X_{\overline{k}_\frakl}/\overline{k}_\frakl}) \!\otimes_{\overline{k}_\frakl}\! \bbC_l(i\!-\!1) \\
\cong& \,\bigoplus_{i=0}^2  H^i(X_{\overline{k}},\Omega^{2-i}_{X_{\overline{k}}/\overline{k}}) \!\otimes_{\overline{k}} \bbC_l(i\!-\!1) \,,
\end{align*}
which is
\mbox{$\bbC_l$-linear}
and equivariant with respect to the natural actions of the decomposition group
$\smash{\Gal(\overline{k}_\frakl/k_\frakl) \subset \Gal(\overline{k}/k)}$.
Restricting to the transcendental part, one obtains a 
$\smash{\Gal(\overline{k}_\frakl/k_\frakl)}$-equivariant
\mbox{$\bbC_l$-linear}
isomorphism
$$\iota_{X,\frakl}^T\colon T_{\overline{l}} \!\otimes_{\overline\bbQ_l}\! \bbC_l \stackrel{\cong}{\longrightarrow} H^0(X_{\overline{k}},\Omega^2_{X_{\overline{k}}/\overline{k}}) \!\otimes_{\overline{k}} \bbC_l(-1) \oplus V_{\overline{k}} \!\otimes_{\overline{k}} \bbC_l \oplus H^2(X_{\overline{k}},\calO_{X_{\overline{k}}}) \!\otimes_{\overline{k}} \bbC_l(1) \,,$$
for~$\smash{V_{\overline{k}} := (c_1(\Pic X_{\overline{k}}) \!\otimes_\bbZ\! \overline{k})^\perp \subset H^1(X_{\overline{k}},\Omega^1_{X_{\overline{k}}/\overline{k}})}$.
\item
When~combining further with the inverse map of the natural isomorphism
$\smash{H^2_\et(X_\bbC, \bbC_l(1)) \stackrel{\cong}{\leftarrow} H^2_\et(X_{\overline{k}}, \bbC_l(1))}$
\cite[Expos\'e~XVI, Corollaire~1.6]{SGA4},
$\iota_{X,\frakl}$
induces a natural~isomorphism
\begin{align*}
H^2_\et(X_\bbC, \bbC_l(1)) \stackrel{\cong}{\longrightarrow}& \,\bigoplus_{i=0}^2 H^i(X_{\overline{k}},\Omega^{2-i}_{X_{\overline{k}}/\overline{k}}) \!\otimes_{\overline{k}}\! \bbC_l(i\!-\!1) \\
\cong& \,\bigoplus_{i=0}^2 H^i(X_\bbC,\Omega^{2-i}_{X_\bbC/\bbC}) \!\otimes_\bbC\! \bbC_l(i\!-\!1) \,.
\end{align*}
Applying the comparison isomorphism
$H^2(X(\bbC), \bbC_l(1)) \stackrel{\cong}{\leftarrow} H^2_\et(X_\bbC, \bbC_l(1))$
\cite[Expos\'e~XI, Th\'eor\`eme~4.4.iii)]{SGA4} to the left and the GAGA isomorphism~\cite[Th\'e\-o\-r\`eme~1]{Se56} to the right hand side, this goes over into the usual Hodge decomposition of complex~cohomology,
\begin{equation}
\label{Hodge}
H^2(X(\bbC), \bbC_l(1)) = \bigoplus_{i=0}^2 H^i(X(\bbC), \Omega^{2-i}_{\an, X(\bbC)}) \!\otimes_\bbC\! \bbC_l(i\!-\!1) \,.
\end{equation}
Note~that the Tate twists do not carry any information here, as
$\bbC_l(i) \cong \bbC_l$,
for every
$i\in\bbZ$,
and there is no Galois action on complex~cohomology.
\end{abc}
\end{ttt}

\subsubsection*{An elementary descent argument}

\begin{lem}[Descent for linear maps]
\label{el_dec}
Let\/~$K$
be a field of
characteristic\/~$0$
and\/
$L \supseteq K$
an algebraically closed~field. For\/
\mbox{$K$-vector}
spaces\/
$V$
and\/~$W$,
let an\/
\mbox{$L$-linear}
map\/~$F\colon V \!\otimes_K\! L \to W \!\otimes_K\! L$
be given that commutes with the actions
of\/~$\Aut_K(L)$.
I.e.~such that, for every
$\sigma \in \Aut_K(L)$,
the diagram
$$
\xymatrix{
V \!\otimes_K\! L \ar@{->}[rr]^F \ar@{->}[d]_\sigma && V \!\otimes_K\! L \ar@{->}[d]^\sigma \\
V \!\otimes_K\! L \ar@{->}[rr]^F && V \!\otimes_K\! L
}
$$
commutes. Then~there is a unique\/
\mbox{$K$-linear}
map\/~$f\colon V \to W$
such that
$F = f \!\otimes_K\! L$.\medskip

\noindent
{\bf Proof.}
{\em
According to \cite[Corollary~53.b)]{Cl}, one has that
$L^{\Aut_K(L)}=K$.
Consequently,~the only
\mbox{$\Aut_K(L)$-invariant}
elements in
$V \!\otimes_K\! L$
are those
in~$V$,
and analogously
for~$W \!\otimes_K\! L$.
By~assumption, one has that
$F(V)$
is
\mbox{$\Aut_K(L)$-invariant},
and therefore
$F(V) \subseteq W$.
The restriction
$f := F |_V \colon V \to W$
is a
\mbox{$K$-linear}
map, as
$V$~is
a
\mbox{$K$-vector}
space and
$F$~is
\mbox{$L$-linear}.
Finally,~it is clear that
$F = f \!\otimes_K\! L$.
}
\eop
\end{lem}

\begin{rem}
This~is certainly a particular case of faithfully flat descent \cite[Expos\'e~VIII, Lemme~1.4]{SGA1}, but Lemma~\ref{el_dec} suffices for the purposes of this~article.
\end{rem}

\section{The various actions of the endomorphism field}

\begin{ttt}[Complex cohomology]
\label{compl_coh}
Let~$X$
be a
$K3$~surface
over a
field~$K$
that is embeddable
into~$\bbC$.
Consider~the embedding
$K \hookrightarrow \bbC$
as being~fixed. The~endomorphism field
$E$
of~$X$
is then defined as~follows.

One~considers the complex manifold
$X(\bbC) = X_\bbC(\bbC)$.
Then~the transcendental part
$T := (H_\alg)^\perp \subset H^2(X(\bbC), \bbQ)$
of the cohomology is a pure
\mbox{weight-$2$}
\mbox{$\bbQ$-Hodge}
structure~\cite[Paragraph~2.1.12 and D\'efinition~2.1.10]{De71}. One
puts~$E := \End_\Hg(T)$
to be the endomorphism ring
of~$T$
in the category of Hodge~structures. It~is known that
$E$
is always a number field~\cite[Theorem~1.6.a)]{Za}, in fact either totally real or a CM~field (cf.\ Paragraph~\ref{CM_defs}, below), at least in the realm of
$K3$~surfaces.
\end{ttt}

\begin{ttt}[Algebraic de Rham cohomology]
\label{AlgdR}
Let~$X$
be as in~\ref{compl_coh}, and let the endomorphism
field~$E$
of~$X$
be as defined~above. Then,~in particular, one has a
\mbox{$\bbQ$-linear}
action
of~$E$
on
$T \subset H^2(X(\bbC), \bbQ)$.
Hence,~$E$
acts
\mbox{$\bbC$-linearly}~on
$$T \!\otimes_\bbQ\! \bbC \subset H^2(X(\bbC), \bbC) = H^2(X(\bbC), \bbQ) \!\otimes_\bbQ\! \bbC \,,$$
and the action commutes with that
of~$\Aut(\bbC)$.
But
\begin{align}
\label{Hodge_dec}
H^2(X(\bbC), \bbQ) \!\otimes_\bbQ\! \bbC = \bigoplus_{i=0}^2 H^i(X(\bbC), \Omega^{2-i}_{\an, X(\bbC)}) \cong \bigoplus_{i=0}^2 &\,H^i(X_\bbC, \Omega^{2-i}_{X/\bbC}) \\[-5mm]
 &\cong \bigoplus_{i=0}^2 H^i(X_{\overline{K}}, \Omega^{2-i}_{X/\overline{K}}) \!\otimes_{\overline{K}}\! \bbC \,, \nonumber
\end{align}
the actions at least
of~$\Aut_{\overline{K}}(\bbC)$
agreeing with each other, as both take place only via the right factor.
The~isomorphism (\ref{Hodge_dec}) lets the
\mbox{$\bbC$-linear}
action
of~$E$
carry over to the transcendental part within the right hand~side. There,~it again commutes with the action
of~$\Aut_{\overline{K}}(\bbC)$.

Thus, Lemma~\ref{el_dec} yields a
\mbox{$\smash{\overline{K}}$-linear}
action
of~$E$
on
$$H^i(X_{\overline{K}}, \Omega^{2-i}_{X/\overline{K}}) \,,$$
for
$i = 0,2$,
as well as on
$$V_{\overline{K}} := (c_1(\Pic X_{\overline\bbC}) \!\otimes_\bbZ\! \overline{K})^\perp = (c_1(\Pic X_{\overline{K}}) \!\otimes_\bbZ\! \overline{K})^\perp \subset H^1(X_{\overline{K}}, \Omega^1_{X/{\overline{K}}}) \,.$$
\end{ttt}

\begin{ttt}[$\overline{l}$-adic
cohomology]
\label{l_adic}
Let~again
$X$
be as in~\ref{compl_coh}, and let the endomorphism
field~$E$
of~$X$
be as defined~above.
Then~for any prime
number~$l$,
the
\mbox{$\bbQ$-linear}
action
of~$E$
on
$T \subset H^2(X(\bbC), \bbQ)$
induces a
\mbox{$\overline\bbQ_l$-linear}
action
of~$E$
on
\begin{align*}
T_{\overline{l}} = T \!\otimes_\bbQ\! \overline\bbQ_l(1) \subset H^2(X(\bbC), \bbQ) \!\otimes_\bbQ\! \overline\bbQ_l(1) &\cong H^2(X(\bbC), \overline\bbQ_l(1)) \\[-1.5mm]
&\;\;\cong H^2_\et(X_\bbC, \overline\bbQ_l(1)) \cong H^2_\et(X_{\overline{K}}, \overline\bbQ_l(1)) \,.
\end{align*}
Note~that the two isomorphisms to the right are canonical
\cite[Expos\'e~XI, Th\'e\-o\-r\`eme~4.4.iii) and Expos\'e~XVI, Corollaire~1.6]{SGA4}.
\end{ttt}

\begin{prop}
\label{transf}
Let\/~$X$
be a\/
$K3$~surface
over a number
field\/~$k$.
Fix~a prime
number\/~$l$
and a
prime\/~$\frakl$
of\/~$k$
lying
above\/~$l$.
Then~the endomorphism
field\/~$E$
of\/~$X$
acts\/
\mbox{$\smash{\overline\bbQ_l}$-linearly}
on\/~$\smash{T_{\overline{l}} \subset H^2_\et(X_{\overline{k}}, \overline\bbQ_l(1))}$
and in such a way that the\/
\mbox{$\smash{\Gal(\overline{k}_\frakl/k_\frakl)}$-equivariant\/}
\mbox{$\bbC_l$-linear}
isomorphism
$$
\iota_{X,\frakl}^T\colon T_{\overline{l}} \!\otimes_{\overline\bbQ_l}\! \bbC_l \stackrel{\cong}{\longrightarrow} H^0(X_{\overline{k}},\Omega^2_{X_{\overline{k}}/\overline{k}}) \!\otimes_{\overline{k}}\! \bbC_l(-1) \oplus V_{\overline{k}} \!\otimes_{\overline{k}}\! \bbC_l \oplus H^2(X_{\overline{k}},\calO_{X_{\overline{k}}}) \!\otimes_{\overline{k}}\! \bbC_l(1)
$$
commutes with the actions
of~$E$
on either~side.\medskip

\noindent
{\bf Proof.}
{\em
The
isomorphism~$\iota_{X,\frakl}^T$
is just the restriction of the isomorphism
$$\iota_{X,\frakl}\colon H^2_\et(X_{\overline{k}}, \bbC_l(1)) \stackrel{\cong}{\longrightarrow} \bigoplus_{i=0}^2  H^i(X_{\overline{k}},\Omega^{2-i}_{X_{\overline{k}}/\overline{k}}) \!\otimes_{\overline{k}} \bbC_l(i\!-\!1) \,,$$
which, according to Paragraph~\ref{nf_case}.b), may be obtained as~follows. Start with the Hodge decomposition~(\ref{Hodge}) of
$H^2(X(\bbC), \bbC_l(1))$
and apply the GAGA isomorphism to the right hand~side. Moreover,~to the left hand side, apply the canonical comparison isomorphism
$\smash{H^2_\et(X_{\overline{k}}, \bbC_l(1)) \stackrel{\cong}{\to} H^2_\et(X_\bbC, \bbC_l(1)) \stackrel{\cong}{\to} H^2(X(\bbC), \bbC_l(1))}$
\cite[Expos\'e~XVI, Corollaire~1.6 and Expos\'e~XI, Th\'eor\`eme~4.4.iii)]{SGA4}.

The~actions
of~$E$
on both sides are constructed from one and the same action on
$H^2(X(\bbC), \bbC_l(1))$
via transport of structure along exactly these isomorphisms. Therefore, they must~agree.
}
\eop
\end{prop}

\section{The actions of the endomorphism field and the Galois group \\ on algebraic de Rham and
$l$-adic
cohomologies---The main results}

\subsubsection*{The general setting--Non-commuting actions.}

\begin{conv}
\label{normiert}
Let~$X$
be a
$K3$~surface
over a
subfield~$K$
of~$\bbC$.
In what follows, we treat the endomorphism field
$E$
as being embedded as a subfield into
$\smash{\overline{K}}$
in the following~way. The~vector space\vspace{.15mm}
$\smash{H^0(X_{\overline{K}},\Omega^2_{X_{\overline{K}}/\overline{K}})}$
is of dimension one, hence the action
of~$E$
is necessarily given~by
${}^{[u]} \omega = \iota(u) \cdot \omega$,
for an embedding
$\smash{\iota\colon E \hookrightarrow \overline{K}}$
and arbitrary
$u \in E$.
We~treat
$E$
as being a subfield
of~$\smash{\overline{K}}$,
via~$\iota$.
Then
$${}^{[u]} \omega = u \cdot \omega \,,$$
for every
$\smash{\omega \in H^0(X_{\overline{K}},\Omega^2_{X_{\overline{K}}/\overline{K}})}$.
\end{conv}

\begin{lem}[Non-commuting actions on algebraic de Rham cohomology]
\label{actions_dR}
Let\/~$X$
be a\/
$K3$~surface
over a
subfield\/
$K$
of\/~$\bbC$
with endomorphism
field\/~$E$.

\begin{abc}
\item
Let\/~$\smash{\omega \in H^0(X_{\overline{K}},\Omega^2_{X_{\overline{K}}/\overline{K}})}$.
Then,~for every\/
$\sigma \in \Gal(\overline{K}/K) \subseteq \Aut_\bbQ(\overline{K})$
and each\/
$u \in E$,
one has\/
$${}^{\sigma \circ [u] \circ \sigma^{-1}} \omega = {}^{[\sigma^{-1}(u)]} \omega \,.$$
\item
Suppose,~in particular, that\/
$\sigma \in \Gal(\overline{K}/K) \subseteq \Aut_\bbQ(\overline{K})$
and\/
$u \in E$
are of the kind that\/
$\sigma(u) \neq u$.
Then~the actions of\/
$\sigma$
and\/~$[u]$
on\/
$\smash{H^0(X_{\overline{K}},\Omega^2_{X_{\overline{K}}/\overline{K}})}$
do not commute with each~other.
\end{abc}\smallskip

\noindent
{\bf Proof.}
{\em
a)
One~has
${}^{[u]} \omega = u \cdot \omega$,
for every
$u \in E$~and
$\smash{\omega \in H^0(X_{\overline{K}},\Omega^2_{X_{\overline{K}}/\overline{K}})}$.
On~the other hand, the action of
$\sigma$
is provided by pull-back,
$${}^\sigma \omega = \sigma^* \omega \,.$$
Consequently,~one~finds
$$
{}^{\sigma \circ [u] \circ \sigma^{-1}} \omega = {}^{\sigma} ({}^{[u]} ({}^{\sigma^{-1}} \omega)) = \sigma^* (u \cdot (\sigma^{-1})^* \omega) = \sigma^{-1}(u) \cdot \sigma^* ((\sigma^{-1})^* \omega) = \sigma^{-1}(u) \cdot \omega \,,
$$
as~claimed. Recall~here that
$\smash{\sigma\colon X_{\overline{K}} \to X_{\overline{K}}}$
is induced by the field automorphism
$\sigma^{-1}\colon \overline{K} \to \overline{K}$.\smallskip

\noindent
b)
As~the assumption implies that
$u \neq \sigma^{-1}(u)$,
part~a) shows that the actions of
\mbox{$\sigma \!\circ\! [u] \!\circ\! \sigma^{-1}$}
and~$[u]$
on~$\smash{H^0(X_{\overline{K}},\Omega^2_{X_{\overline{K}}/\overline{K}})}$
do not~coincide. This~is equivalent to the assertion.
}
\eop
\end{lem}

\subsubsection*{A lower bound for the component group of the algebraic monodromy group}

\begin{prop}
\label{nicht_neutral}
Let\/
$X$
be a\/
$K3$~surface
over a number
field\/~$k$
with endomorphism
field\/~$E$.
Moreover,~let\/
$\sigma \in \Gal(\overline{k}/k) \!\setminus\! \Gal(\overline{k}/kE)$
be~arbitrary. Then one~has
$$\varrho_{X,\overline{l}}(\sigma) \not\in G_{X,\overline{l}}^0 \,,$$
for any prime
number\/~$l$.\medskip

\noindent
{\bf Proof.}
{\em
{\em First step.}
Adjusting the
automorphism~$\sigma$.
Fixing~a prime
number~$l$.\smallskip

\noindent
Let~$k' \supseteq k$
denote the splitting field
of~$\smash{G_{X,\overline{l}}/G_{X,\overline{l}}^0}$,
which is a normal extension and, according to Corollary~\ref{split_l}, independent
of~$l$.
Moreover,~consider the normal closure
$\smash{(kE)^{(n)}}$
of~$kE$
as an extension field
of~$k$.
Then~$\smash{(kE)^{(n)} \supseteq k}$
and hence
$\smash{k'(kE)^{(n)} \supseteq k}$,
too, are normal extensions. It~follows that the assumption, as well as the assertion, depend only on the restriction
$\smash{\sigma |_{k'(kE)^{(n)}} \in \Gal(k'(kE)^{(n)}/k)}$.

Since~$\smash{k'(kE)^{(n)} \supseteq k}$
is an extension of number fields, the Chebotarev density theorem provides a prime
$\frakL$
of~$\smash{k'(kE)^{(n)}}$
that is unramified
over~$k$
such that
$\smash{\Frob^\frakL_{k'(kE)^{(n)}/k} = \sigma |_{k'(kE)^{(n)}}}$.
We~denote the prime
of~$k$
lying under
$\frakL$
by~$\frakl$.
Let~us choose a further extension
$\frakM$
of~$\frakL$
to~$\smash{\overline{k}}$
and put
$$\widetilde\sigma := \Frob^\frakM_{\overline{k}/k} \in \Gal(\overline{k_\frakl}/k_\frakl) \hookrightarrow \Gal(\overline{k}/k) \,.$$
Then~$\smash{\widetilde\sigma |_{k'(kE)^{(n)}} = \sigma |_{k'(kE)^{(n)}} \in \Gal(k'(kE)^{(n)}/k)}$,
so that it suffices to show the assertion
for~$\smash{\widetilde\sigma}$,
instead
of~$\sigma$.
Note~at this point that the choice
of~$\frakM$
distinguishes a particular embedding
$\smash{\Gal(\overline{k_\frakl}/k_\frakl) \hookrightarrow \Gal(\overline{k}/k)}$.

Let~$l$
be the prime number lying under the
prime~$\frakl$
of~$k$.
In~the steps below,
\mbox{$l$-adic}
cohomology is always meant to be for
this value
of~$l$.\medskip

\noindent
{\em Second step.}
The actions on algebraic de Rham cohomology.\smallskip

\noindent
By~assumption, there is an element
$u_0 \in E$
such that
$\widetilde\sigma(u_0) \neq u_0$.
Then Lemma \ref{actions_dR}.b) yields that the actions of
$\widetilde\sigma$
and~$[u_0]$
on
$$\smash{H^0(X_{{\overline{k}_\frakl}},\Omega^2_{X_{{\overline{k}_\frakl}}/\overline{k}_\frakl})}$$
do not commute with each~other. Consequently,~the actions of
\mbox{$\smash{\widetilde\sigma \!\circ\! [u_0] \!\circ\! \widetilde\sigma^{-1}}$}
and
$[u_0]$~on
$$H^0(X_{{\overline{k}_\frakl}},\Omega^2_{X_{{\overline{k}_\frakl}}/\overline{k}_\frakl}) \!\otimes_{\overline{k}_\frakl}\! \bbC_l(-1) \oplus V_{{\overline{k}_\frakl}} \!\otimes_{\overline{k}_\frakl}\! \bbC_l \oplus H^2(X_{{\overline{k}_\frakl}},\calO_{X_{{\overline{k}_\frakl}}}) \!\otimes_{\overline{k}_\frakl}\! \bbC_l(1) \,,$$
which are both
\mbox{$\bbC_l$-linear}
maps, must be different, as~well.\medskip

\noindent
{\em Third step.}
The transfer to
\mbox{$l$-adic}
cohomology.\smallskip

\noindent
According~to Proposition~\ref{transf}, the actions of
$\widetilde\sigma \!\circ\! [u_0] \!\circ\! \widetilde\sigma^{-1}$
and~$[u_0]$~on
$$T_{\overline{l}} \!\otimes_{\overline\bbQ_l}\! \bbC_l \,,$$
which are again both
\mbox{$\bbC_l$-linear}
maps, do not coincide~either. This~implies that the underlying
\mbox{$\overline\bbQ_l$-linear}
endomorphisms
of~$T_{\overline{l}}$
already differ from each~other. In~other words, the action
of~$\smash{\widetilde\sigma \in \Gal(\overline{k_\frakl}/k_\frakl)}$
on~$T_{\overline{l}}$
does not commute with that
of~$u_0 \in E$.

I.e.,
$\varrho_{X,\overline{l}}(\widetilde\sigma) \not\in \C_E(\O(T_{\overline{l}}))$.
In~particular, one certainly has
$\varrho_{X,\overline{l}}(\widetilde\sigma) \not\in (\C_E(\O(T_{\overline{l}})))^0$,
which, by the work of S.\,G.\ Tankeev~\cite{Tan90,Tan95} and Yu.\,G.\ Zarhin~\cite[Theorem 2.2.1]{Za}, is equivalent to
$\smash{\varrho_{X,\overline{l}}(\widetilde\sigma) \not\in G_{X,\overline{l}}^0}$.
Cf.~formula~(\ref{Tankeev}) from the~introduction.\medskip

\noindent
{\em Fourth step.}
Conclusion.\smallskip

\noindent
Let~$l' \neq l$
be an arbitrary prime~number. Then~one has
$\smash{\varrho_{X,\overline{l}'}(\widetilde\sigma) \not\in G_{X,\overline{l}'}^0}$,
too, by Theorem~\ref{Serre_Lettre}, which completes the~proof.
}
\eop
\end{prop}

\begin{coro}
\label{split_E}
Let\/
$X$
be a\/
$K3$~surface
over a number
field\/~$k$
with endomorphism
field\/~$E$.
Then,~for any prime
number\/~$l$,
the splitting field
of\/~$\smash{G_{X,\overline{l}}/G_{X,\overline{l}}^0}$
contains\/~$kE$
as a~subfield.
\eop
\end{coro}

\begin{theo}
\label{main}
Let\/
$X$
be a\/
$K3$~surface
over a number
field\/~$k$
with endomorphism
field\/~$E$
and
let\/~$l$
be any prime~number.

\begin{abc}
\item
There is a natural surjective homomorphism\/
$\smash{G_{X,\overline{l}}/G_{X,\overline{l}}^0 \twoheadrightarrow \Gal((kE)^{(n)}/k)}$.
If\/~$kE \supseteq k$
is a normal extension then one has a natural surjective homomorphism
\begin{equation}
\label{hom_main_Haupttext}
G_{X,\overline{l}}/G_{X,\overline{l}}^0 \twoheadrightarrow \Gal(kE/k) \,.
\end{equation}
\item
In particular,
$\smash{\#(G_{X,\overline{l}}/G_{X,\overline{l}}^0)}$
is always a multiple
of\/~$[kE:k]$.
\end{abc}\smallskip

\noindent
{\bf Proof.}
{\em
a)
Write
$k' \supseteq k$
for the splitting field
of~$\smash{G_{X,\overline{l}}/G_{X,\overline{l}}^0}$.
Then~$k'$
is normal
over~$k$
and one has a natural~isomorphism
$$\smash{G_{X,\overline{l}}/G_{X,\overline{l}}^0 \stackrel{\cong}{\longrightarrow} \Gal(k'/k) \,.}$$
Corollary~\ref{split_E} implies that
$k' \supseteq kE$.
Therefore,~it follows that
$\smash{k' \supseteq (kE)^{(n)}}$
holds, too, which implies the first~claim. The~second claim follows~immediately.\smallskip

\noindent
b)
is a direct consequence of~a).
}
\eop
\end{theo}

\begin{rem}
At~least under the assumption of the Hodge conjecture,
$kE \supseteq k$
is always a normal extension of~fields. Cf.~Theorem~\ref{normality}.
\end{rem}

\section{The fundamental lemma on the CM case}

\begin{ttt}
\label{CM_defs}
Recall~that a number
field~$E$
is said to be a {\em CM field\/} if
$E$~is
a totally imaginary quadratic extension of a
totally real~field.
Let~$E$
be a CM field and let
$u_0 \in E$
be a primitive~element. Then
$\Gal(E^{(n)}/\bbQ)$
acts transitively on the conjugates
$u_1, \overline{u}_1, \ldots, u_d, \overline{u}_d$
of~$u_0$,
thereby preserving the obvious system of blocks of
size~$2$.
Moreover,~the complex conjugation just flips each block, and is therefore a central element
in~$\Gal(E^{(n)}/\bbQ)$
\cite[Proposition~5.11]{Sh}.

A~$K3$~surface
$X$
over a number
field~$k$
is said to have {\em complex multiplication (CM)\/} if the endomorphism field
$E$
is a CM~field. In~this situation, let
$E_0 \subset E$
be the maximal totally real subfield,
$d := [E_0:\bbQ]$,
and~$\smash{r := \dim_{\overline\bbQ} T}$.
Then,~for every prime
number~$l$,
the action
of~$u_0$
splits~$T_{\overline{l}}$
into eigenspaces
$T_{\overline{l}} = T_{l,u_1} \oplus T_{l,\overline{u}_1} \oplus \cdots \oplus T_{l,u_d} \oplus T_{l,\overline{u}_d}$,
which come in pairs corresponding to the conjugates
of~$u_0$.
Each~eigenspace is of
dimension~$\frac{r}{2d}$.
Indeed,~in
$T \subset H^2(X(\bbC), \bbQ)$,
the eigenspaces are mutually conjugate, and hence of the same dimension. Thus,~the~claim follows from Proposition~\ref{transf}.
%
\end{ttt}

\begin{ttt}
\label{jump_char}
Let,~as before,
$X$
be a
$K3$~surface
over a number
field~$k$.
Then~the {\em jump character\/}
of~$X$
(cf.\ \cite[Definition~2.4.6]{CEJ}) is the homomorphism
$\smash{\tau_X\colon \Gal(\overline{k}/k) \to \{\pm1\}}$
given by the natural action
of~$\smash{\Gal(\overline{k}/k)}$
on
$\Exterior^\maxi T_{\overline{l}}$.
\end{ttt}

\begin{rems}
\begin{iii}
\item
The~jump character is independent
of~$l$
\cite[Proposition 2.2.1.b)]{CEJ}.
\item
The name ``jump character'' has not been chosen at random. In fact, suppose that
$\rk\Pic X_{\overline{k}}$
is even and let
$\frakp$
be a prime of
$k$,
at which
$X$
has good reduction and such that the jump character
of~$X$
evaluates
to~$(-1)$
at~$\Frob_\frakp$.
Then
$\smash{\rk\Pic X_{\overline\bbF_{\!\frakp}} \geq \rk\Pic X_{\overline{k}} + 2}$
\cite[Proposition~2.4.2]{CEJ}.
\end{iii}
\end{rems}

\begin{ttt}
Contrary~to the generic and RM cases, in the CM~case, one has that (\ref{hom_main_Haupttext}) is always an~isomorphism. I.e., that the splitting field
of~$\smash{G_{X,\overline{l}}/G_{X,\overline{l}}^0}$
is
exactly~$kE$.
Cf.~Theorem~\ref{jump_CM_alg}.a). As~a consequence, it is possible in the CM case to describe the jump character entirely in terms of the endomorphism field and the Picard rank, cf.~Theorem~\ref{jump_CM_alg}.b).

The~goal of this section is to provide the necessary preparations. The complete argument will be given in Section~\ref{Hodge_dep} and is, unfortunately, somewhat indirect. It~makes use of the Hodge~conjecture, which is known to be true in the case of CM
$K3$~surfaces.
\end{ttt}

\begin{prop}
\label{CM_general}
Let\/~$X$
be a\/
$K3$~surface
over a number
field\/~$k$
having CM by an endomorphism
field\/~$E$.
Choose~a prime
number\/~$l$.

\begin{abc}
\item
Then~each of the eigenspaces\/
$T_{l,u_1}, T_{l,\overline{u}_1}, \ldots, T_{l,u_d}, T_{l,\overline{u}_d}$
is~isotropic. One~has\/
$\C_E(\O(T_{\overline{l}})) \cong [\GL_{r/2d}(\overline\bbQ_l)]^d$.
In~particular,
$G_{X,\overline{l}}^0 \cong \C_E(\O(T_{\overline{l}}))$.
\item
There~is a natural faithful permutation action
of\/~$\smash{G_{X,\overline{l}}/G_{X,\overline{l}}^0}$
on the set of eigen\-spaces\/
$\{T_{l,u_1}, T_{l,\overline{u}_1}, \ldots, T_{l,u_d}, T_{l,\overline{u}_d}\}$
that preserves the obvious system of blocks of
size\/~$2$,
$$\smash{G_{X,\overline{l}}/G_{X,\overline{l}}^0 \hookrightarrow (\bbZ/2\bbZ)^d \!\rtimes\! S_d \subseteq S_{2d} \,.}$$
\end{abc}

\noindent
{\bf Proof.}
{\em
a)
For~two eigenvectors,
$\langle v,v' \rangle \neq 0$
is possible only when the corresponding eigenvalues
$u,u'$
are complex~conjugates. Indeed one~has
$$\smash{u \!\cdot\! \langle v,v' \rangle  = \langle {}^{[u_0]} v, v' \rangle  = \langle v, {}^{[\overline{u}_0]} v' \rangle  = \overline{u}' \!\cdot\! \langle v,v' \rangle \,,}$$
due to \cite[Theorem~1.5.1]{Za}. In~particular, the first claim is~shown.

Moreover, the cup product pairing
on~$T_{\overline{l}}$
is nondegenerate, so there exist bases
$B_{u_1}, B_{\overline{u}_1}, \ldots, B_{u_d}, B_{\overline{u}_d}$
of 
$T_{l,u_1}$,
$T_{l,\overline{u}_1}$,
\ldots,
$T_{l,u_d}$,
and~$T_{l,\overline{u}_d}$,
respectively, such that,
for~$B := B_{u_1} \!\cup\! B_{\overline{u}_1} \!\cup\cdots\cup\! B_{u_d} \!\cup\! B_{\overline{u}_d}$,
the matrix
$\M_B(\langle .,. \rangle)$
is block diagonal consisting
of~$d$
diagonal blocks of~type
$$\textstyle (\atop{\;\;\;0 \;\;\;E_{r/2d}}{E_{r/2d} \;\;\;0\;\;\;})$$
and zeroes,~otherwise.

Furthermore,~an
endomorphism~$s \in \End(T_{\overline{l}})$
commutes with the action
of~$E$
if and only if it maps the eigenspaces to~themselves. I.e., when
{\arraycolsep0pt
$$\M_B^B(s) =
\left(
\begin{array}{*{5}{>{\scriptscriptstyle}c}}
G_1^+  & 0      & \ldots & 0       & 0      \\[-1.3mm]
0      & G_1^-  & \ldots & 0       & 0      \\[-1.3mm]
\ldots & \ldots & \ldots & \ldots  & \ldots \\[-1.3mm]
0      & 0      & \ldots & \;G_d^+ & 0      \\[-1.3mm]
0      & 0      & \ldots & 0       & \!G_d^-
\end{array}
\right)
\,,$$}%
for suitable matrices
$G_1^+, G_1^-, \ldots, G_d^+, G_d^-$.
On~the other hand, a direct calculation shows that
$s$~is
orthogonal if and only if
$\smash{G_i^- = ((G_i^+)^t)^{-1}}$,
for
each~$i$.
Consequently,
{\arraycolsep0pt
\begin{align*}
\C_E(\O(T_{\overline{l}})) &= \left\{s \in \End(T_{\overline{l}}) \,\Biggg|\, \M_B^B(s) =
\left(
\begin{array}{*{5}{>{\scriptscriptstyle}c}}
G_1    & 0                & \ldots & \;\;0      & 0               \\[-1.3mm]
0      & \!(G_1^t)^{-1}\! & \ldots & \;\;0      & 0               \\[-1.3mm]
\ldots & \ldots           & \ldots & \;\;\ldots & \ldots          \\[-1.3mm]
0      & 0                & \ldots & \;\;\;G_d  & 0               \\[-1.3mm]
0      & 0                & \ldots & \;\;0      & \!\!(G_d^t)^{-1}
\end{array}
\right)
{\rm ,\,for~some~} G_1, \ldots, G_d \right\} \\
 &\cong [\GL_{r/2d}(\overline\bbQ_l)]^d \,,
\end{align*}}%
as~claimed.
In~particular,
$\C_E(\O(T_{\overline{l}}))$
is connected, which implies the final assertion in view of~(\ref{Tankeev}).\smallskip

\noindent
b)
As~the neutral component of an algebraic group is always normal, one clearly has
$$\smash{G_{X,\overline{l}} \subseteq \N_{\O(T_{\overline{l}})} (G_{X,\overline{l}}^0) \cong \N_{\O(T_{\overline{l}})} (\C_E(\O(T_{\overline{l}}))) \,.}$$
Moreover,~if\vspace{.2mm}
$s\in\O(T_{\overline{l}})\subset\End(T_{\overline{l}})$
normalises
$\smash{G_{X,\overline{l}}^0 \cong \C_E(\O(T_{\overline{l}}))}$
then each of the subvector spaces
$\smash{s(T_{l,u_1}), s(T_{l,\overline{u}_1}), \ldots, s(T_{l,u_d}), s(T_{l,\overline{u}_d})}$
must be invariant under the action
of~$\smash{\C_E(\O(T_{\overline{l}}))}$.
Since~there are no such spaces other than
$\smash{T_{l,u_1}, T_{l,\overline{u}_1}, \ldots, T_{l,u_d}, T_{l,\overline{u}_d}}$,
the endomorphism
$s$
necessarily permutes them, thereby defining a permutation action
\begin{equation}
\label{perm_act}
G_{X,\overline{l}} \longrightarrow \Sym(\{T_{l,u_1}, T_{l,\overline{u}_1}, \ldots, T_{l,u_d}, T_{l,\overline{u}_d}\}) \cong S_{2d} \,.
\end{equation}
Furthermore,~$s\in\O(T_{\overline{l}})\subset\End(T_{\overline{l}})$
commutes
with~$E$
if and only if
$s$
maps every eigenspace to itself, so the kernel of~(\ref{perm_act}) is exactly
$\smash{\C_E(\O(T_{\overline{l}})) \cong G_{X,\overline{l}}^0}$,
which implies the first~assertion. The preservation of the block system is~obvious.
}
\eop
\end{prop}

\begin{lem}[The fundamental lemma on the CM~case]
\label{fund_lem}
Let\/~$X$
be a\/
$K3$~surface
over a number
field\/~$k$
having CM by an endomorphism
field\/~$E$.
For~some prime
number\/~$l$,
suppose~that
\begin{equation}
\label{non_comm}
{}^{\sigma \circ [u] \circ \sigma^{-1}} v = {}^{[\sigma^{-1}(u)]} v \,,
\end{equation}
for every\/
$v \in T_{\overline{l}}$,
$\smash{\sigma \in \Gal(\overline{k}/k)}$,
and\/~$u \in E$.

\begin{abc}
\item
Then~there is a natural isomorphism
$$G_{X,\overline{l}}/G_{X,\overline{l}}^0 \stackrel{\cong}{\longrightarrow} \Gal(kE/k) \,.$$
I.e.,~the splitting field
of\/~$\smash{G_{X,\overline{l}}/G_{X,\overline{l}}^0}$
is
exactly\/~$kE$.
\item
The~jump character\/
$\tau_X\colon \Gal(\overline{k}/k) \to \{\pm1\}$
is
$$
\sigma \mapsto \left\{
\begin{array}{rrr}
1 \,,                                                                             & \hspace{.7cm}\textstyle\text{~if\/~} \smash{\frac{22-\rk\Pic X_{\overline{k}}}{[E:\bbQ]}} \text{~is~even} \,, \\
\text{sign~of~the~natural~permutation~}                                           & \\[-.4mm]
\text{action~of~} \sigma \in \smash{\Gal(\overline{k}/k)} \text{~on~the~}         & \\[-.4mm]
\text{conjugates\/~} u_1, \overline{u}_1, \ldots, u_d, \overline{u}_d \text{~of~} & \\[-.4mm]
\text{a~primitive~element~} u_0 \in E \,,                                         & \textstyle\text{~if\/~} \smash{\frac{22-\rk\Pic X_{\overline{k}}}{[E:\bbQ]}} \text{~is~odd} \,. \;\;
\end{array}
\right.
$$
\end{abc}

\noindent
{\bf Proof.}
{\em
a)
The~assertion means that
$\smash{\varrho_{X,\overline{l}}(\sigma) \in G_{X,\overline{l}}^0}$
if and only if
$\smash{\sigma\in\Gal(\overline{k}/kE)}$,
which is equivalent~to
\begin{equation}
\label{comm}
\varrho_{X,\overline{l}}(\sigma) \in G_{X,\overline{l}}^0 \;\Longleftrightarrow\; \sigma |_E = \id_E \,.
\end{equation}
The implication
``$\Longrightarrow$''
of~(\ref{comm}) is true, due to Theorem~\ref{main}.a). In~order to show
``$\Longleftarrow$'',
note that, for
$\sigma |_E = \id_E$,
assumption~(\ref{non_comm}) says that the action
of~$\sigma$
on~$T_{\overline{l}}$
commutes with that
of~$E$.
But~this is sufficient to imply
$\smash{\varrho_{X,\overline{l}}(\sigma) \in G_{X,\overline{l}}^0}$
in the CM case, due to Proposition~\ref{CM_general}.a).
For the statement concerning the splitting field of
$\smash{G_{X,\overline{l}}/G_{X,\overline{l}}^0}$,
recall part~\ref{splitF} of the Conventions and Notation subsection.\smallskip

\noindent
b)
Assumption~(\ref{non_comm}) yields
$\smash{{}^{\sigma \circ [\sigma(u_0)]} v = {}^{[u_0]} ({}^\sigma v)}$,
for every
$\smash{\sigma \in \Gal(\overline{k}/k)}$
and~$v \in T_{\overline{l}}$.
Hence,~for an eigenvector
$\smash{v \in T_{l,\nu}}$,
one has
${\smash{}^\sigma v \in T_{l,\sigma(\nu)}}$.
In~other words, the action
of~$\smash{\Gal(\overline{k}/k)}$
permutes the eigenspaces
$T_{l,u_1}, T_{l,\overline{u}_1}, \ldots, T_{l,u_d}, T_{l,\overline{u}_d}$
according to the natural action
of~$\smash{\Gal(\overline{k}/k)}$
on
$u_1, \overline{u}_1, \ldots, u_d, \overline{u}_d$,
as elements
of~$\smash{\overline{k}}$.

On~the other hand, by definition, the jump character maps each
$\smash{\sigma\in\Gal(\overline{k}/k)}$
to
$\smash{\det\varrho_{X,\overline{l}}(\sigma)}$.
As~$\det$
is locally constant, this depends only on the class
of~$\smash{\varrho_{X,\overline{l}}(\sigma) \in G_{X,\overline{l}}}$
in the factor
group~$\smash{G_{X,\overline{l}}/G_{X,\overline{l}}^0 \hookrightarrow (\bbZ/2\bbZ)^d \!\rtimes\! S_d}$.
By~what was just shown, this class is given directly by the natural action
of~$\sigma$
on~$u_1, \overline{u}_1, \ldots, u_d, \overline{u}_d$.

We~claim, more generally, that, for every
$A \in \N_{\O(T_{\overline{l}})} (\C_E(\O(T_{\overline{l}})))$,
one~has
$$\det A = (\sgn (\pi_A, a_A))^{r/2d} \,,$$
where
$(\pi_A, a_A)$
denotes the class
of~$A$
in
$$\N_{\O(T_{\overline{l}})} (\C_E(\O(T_{\overline{l}})))/\C_E(\O(T_{\overline{l}})) \cong (\bbZ/2\bbZ)^d \!\rtimes\! S_d \subseteq S_{2d} \,.$$
Note~here that
$\frac{r}{2d} = \frac{22-\rk\Pic X_{\overline{k}}}{[E:\bbQ]}$.

For~this, choose bases
$\smash{B_1 = \{v_{1,1},\ldots,v_{1,r/2d}\}}$,
\ldots,
$\smash{B_d = \{v_{d,1},\ldots,v_{d,r/2d}\}}$
of
$\smash{T_{l,u_1}, \ldots, T_{l,u_d}}$,
respectively, and equip
$\smash{T_{l,\overline{u}_i}}$
with the basis
$\smash{B_i^* = \{v_{1,1}^*,\ldots,v_{1,r/2d}^*\}}$,
dual
to~$B_i$,
for~$i = 1,\ldots,d$.
Moreover,~for every
$(\pi,a) \in (\bbZ/2\bbZ)^d \!\rtimes\! S_d \subseteq S_{2d}$,
let~$M_{(\pi,a)}$
be the
\mbox{$\bbQ_l$-linear}
map that sends
$\smash{v_{i,j}}$
to
$\smash{v_{\pi(i),j}}$
or
$\smash{v_{\pi(i),j}^*}$,
depending on whether
$a_i = 0$
or~$1$,
and
$\smash{v_{i,j}^*}$
to
$\smash{v_{\pi(i),j}^*}$
or
$\smash{v_{\pi(i),j}}$,
accordingly.
Then,~by construction,
$M_{(\pi,a)}$
is an orthogonal map, hence
$M_{(\pi,a)} \in \N_{\O(T_{\overline{l}})} (\C_E(\O(T_{\overline{l}})))$.
On~the other hand,
$M_{(\pi,a)}$
is the
\mbox{$\bbQ_l$-linear}
map corresponding to a permutation matrix, for the permutation given by
$\smash{\frac{r}{2d}}$
disjoint copies
of~$(\pi,a)$.
Therefore,~$\smash{\det M_{(\pi,a)} = (\sgn (\pi,a))^{r/2d}}$.
Finally,~for
$\smash{A \in \N_{\O(T_{\overline{l}})} (\C_E(\O(T_{\overline{l}})))}$
arbitrary,
$M_{(\pi_A,a_A)}$
belongs to the same component of the algebraic group
$\smash{\N_{\O(T_{\overline{l}})} (\C_E(\O(T_{\overline{l}})))}$
as~$A$.
This~finally shows that
$\smash{\det A = \det M_{(\pi_A,a_A)} = (\sgn(\pi_A,a_A))^{r/2d}}$,
as~claimed.
}
\eop
\end{lem}

\begin{rem}
Supposedly, assumption~(\ref{non_comm}) is always true, even in the non-CM~cases. Our~argument, however, relies on the Hodge~conjecture and is therefore postponed to Section~\ref{Hodge_dep}. Cf.\ Lemma~\ref{comm-dR}.b), below. 
\end{rem}

\section{A digression: the field of definition of the conjectural correspondence describing real or complex multiplication}

\begin{ttt}
Let~$X$
be a
$K3$~surface
over a
field~$K\subseteq\bbC$.
Then~one has the direct decomposition
$H^2(X(\bbC), \bbQ) = H_\alg \oplus T$
of pure
\mbox{weight-$2$}
\mbox{$\bbQ$-Hodge}
structures.

Moreover,~let
$E$
be the endomorphism field
of~$X$.
Then~every
$u \in E$
defines a morphism
$[u]\colon T \to T$
of Hodge~structures.
There~exists a (non-unique) extension
$H^2(X(\bbC), \bbQ) \to H^2(X(\bbC), \bbQ)$
of~$[u]$
that is still a morphism of Hodge~structures. Let~us fix one such extension, which, by abuse of notation, we again denote
by~$\smash{[u]}$.

Then~$\smash{[u]\colon H^2(X(\bbC), \bbQ) \to H^2(X(\bbC), \bbQ)}$
is a
\mbox{$\bbQ$-linear}
map and the induced
\mbox{$\bbC$-linear}
map
$H^2(X(\bbC), \bbC) \to H^2(X(\bbC), \bbC)$,
which respects the Hodge decomposition
$\smash{H^2(X(\bbC), \bbC) = \bigoplus_{i=0}^2 H^i(X(\bbC), \Omega^{2-i}_{\an, X(\bbC)})}$,
clearly commutes with the action
of~$\Aut_{\overline{K}}(\bbC)$
on both~sides. Therefore, Lemma~\ref{el_dec} implies that
$[u]$
descends from a
\mbox{$\bbC$-linear}
endomorphism of
$\smash{H^i(X(\bbC), \Omega^{2-i}_{\an, X(\bbC)}) = H^i(X_\bbC, \Omega^{2-i}_{X_\bbC/\bbC})}$
to a
\mbox{$\smash{\overline{K}}$-linear}
endomorphism of
$\smash{H^i(X_{\overline{K}}, \Omega^{2-i}_{X_{\overline{K}}/\overline{K}})}$,
for~$i=0,1,2$.
Cf.\ Paragraph~\ref{AlgdR}.
\end{ttt}

\begin{defi}
Let~$X$
be a
$K3$~surface
with endomorphism
field~$E$
over a field
$K\subseteq\bbC$,
and
let~$u \in E$
be~arbitrary.

\begin{abc}
\item
For~a basis
$(v_1^\an, \ldots, v_{22}^\an)$
of~$H^2(X(\bbC), \bbQ)$
and
$((v_1^\an)^*, \ldots, (v_{22}^\an)^*)$
the dual basis, the {\em analytic Casimir element\/}
$\Delta^\an_u$
associated
with~$u$
is given by
\begin{align*}
\Delta^\an_u := (v_1^\an)^* \!\otimes {}^{[u]} v_1^\an + \cdots + (v_{22}^\an)^* \!\otimes {}^{[u]} v_{22}^\an &\in H^2(X(\bbC), \bbQ)^* \otimes_\bbQ H^2(X(\bbC), \bbQ) \\
 &\cong H^2(X(\bbC), \bbQ) \otimes_\bbQ H^2(X(\bbC), \bbQ) \\
 &\subset H^4((X \!\times\! X)(\bbC), \bbQ) \,.
\end{align*}
\item
For~a basis
$(v_1, \ldots, v_{22})$
of~$\smash{\bigoplus_{i=0}^2 \!H^i(X_{\overline{K}}, \Omega^{2-i}_{X_{\overline{K}}/\overline{K}})}$
and
$(v_1^*, \ldots, v_{22}^*)$
the dual basis, the {\em algebraic Casimir element\/}
$\Delta_u$
associated
with~$u$
is given by
\begin{align}
\label{Casimir}
\Delta_u := v_1^* \otimes {}^{[u]} v_1 + \cdots + v_{22}^* \otimes {}^{[u]} v_{22} &\in \bigoplus_{i=0}^2 H^i(X_{\overline{K}}, \Omega^{2-i}_{X_{\overline{K}}/\overline{K}})^* \otimes_{\overline{K}} \bigoplus_{i=0}^2 H^i(X_{\overline{K}}, \Omega^{2-i}_{X_{\overline{K}}/\overline{K}}) \\
 &\cong \bigoplus_{i=0}^2 H^i(X_{\overline{K}}, \Omega^{2-i}_{X_{\overline{K}}/\overline{K}}) \otimes_{\overline{K}} \bigoplus_{i=0}^2 H^i(X_{\overline{K}}, \Omega^{2-i}_{X_{\overline{K}}/\overline{K}}) \nonumber \\
 &\subset \bigoplus_{i=0}^4 H^i((X \!\times\! X)_{\overline{K}}, \Omega^{4-i}_{(X \times X)_{\overline{K}}/\overline{K}}) \,, \nonumber
\end{align}
Here,
$$[u]\colon \bigoplus_{i=0}^2 H^i(X_{\overline{K}}, \Omega^{2-i}_{X_{\overline{K}}/\overline{K}}) \longrightarrow \bigoplus_{i=0}^2 H^i(X_{\overline{K}}, \Omega^{2-i}_{X_{\overline{K}}/\overline{K}}) \,,
\qquad
v \mapsto {}^{[u]} v
$$
is~the endomorphism described~above.
\end{abc}
\end{defi}

\begin{rems}
\begin{iii}
\item
The~Casimir elements
$\Delta_u$
and~$\Delta^\an_u$
do not depend on the bases~chosen. Cf.~\cite[Section~6.2]{Hu}, where this fact is shown in a rather different~context.
\item
In~particular, the image
of~$\Delta^\an_u$~in
$$H^4((X \!\times\! X)(\bbC), \bbC) \cong \bigoplus_{i=0}^4 H^i((X \!\times\! X)_\bbC, \Omega^{4-i}_{(X \times X)_\bbC/\bbC})$$
under change of coefficients may be constructed in the same way from any basis~of
$$H^2(X(\bbC), \bbC) = \bigoplus_{i=0}^2 H^i(X(\bbC), \Omega^{2-i}_{\an, X(\bbC)}) = \bigoplus_{i=0}^2 H^i(X_\bbC, \Omega^{2-i}_{X_\bbC/\bbC}) \,.$$
Therefore,~it coincides with the image
of~$\Delta_u$
in~$\smash{\bigoplus_{i=0}^4 H^i((X \!\times\! X)_\bbC, \Omega^{4-i}_{(X \times X)_\bbC/\bbC})}$
under base~change.
\end{iii}
\end{rems}

\begin{lem}
\label{Casimir_sigma}
Let\/~$X$
be a\/
$K3$~surface
with endomorphism
field\/~$E$
over a field\/
$K\subseteq\bbC$,
and
let\/~$u \in E$
be~arbitrary. Then,~for any\/
$\sigma\in\Gal(\overline{K}/K)$
and any basis\/
$(v_1, \ldots, v_{22})$
of\/~$\smash{\bigoplus_{i=0}^2 H^i(X_{\overline{K}}, \Omega^{2-i}_{X_{\overline{K}}/\overline{K}})}$,
one~has
$${}^\sigma (\Delta_u) = v_1^* \otimes {}^{(\sigma \circ [u] \circ \sigma^{-1})} v_1 + \cdots + v_{22}^* \otimes {}^{(\sigma \circ [u] \circ \sigma^{-1})} v_{22} \,.$$
%

\noindent
{\bf Proof.}
{\em
Recall~that
$\Delta_u$
is independent of the basis~chosen. Instead~of
$(v_1, \ldots, v_{22})$,
another basis is provided by
$\smash{({}^{\sigma^{-1}} v_1, \ldots, {}^{\sigma^{-1}} v_{22})}$.
Moreover,~as
$\sigma^{-1}$
acts by an orthogonal map, the dual basis is then simply
$\smash{({}^{\sigma^{-1}} (v_1^*), \ldots, {}^{\sigma^{-1}} (v_{22}^*))}$.
Consequently, one~has
$$\Delta_u = {}^{\sigma^{-1}} (v_1^*) \otimes {}^{[u]} ({}^{\sigma^{-1}} v_1) + \cdots + {}^{\sigma^{-1}} (v_{22}^*) \otimes {}^{[u]} ({}^{\sigma^{-1}} v_{22}) \,,$$
which yields~that
\begin{align*}
{}^\sigma (\Delta_u) &= v_1^* \otimes {}^\sigma ({}^{[u]} ({}^{\sigma^{-1}} v_1)) + \cdots + v_{22}^* \otimes {}^\sigma ({}^{[u]} ({}^{\sigma^{-1}} v_{22})) \nonumber \\
 &= v_1^* \otimes {}^{(\sigma \circ [u] \circ \sigma^{-1})} v_1 + \cdots + v_{22}^* \otimes {}^{(\sigma \circ [u] \circ \sigma^{-1})} v_{22} \,,
\end{align*}
as~claimed.
}
\eop
\end{lem}

\begin{prop}
\label{corr}
Let\/~$X$
be a\/
$K3$~surface
with endomorphism
field\/~$E$
over a field\/
$K\subseteq\bbC$,
and
let\/~$u \in E$
be~arbitrary. Suppose~that the Hodge conjecture~\cite{De06} holds
for\/~$(X \!\times\! X)(\bbC)$.

\begin{abc}
\item
Then~there exists a\/
\mbox{$\bbQ$-cycle\/}
$C_u$
of
codimension\/~$2$
in\/~$(X \!\times\! X)_{\overline{K}}$
of the kind~that
\begin{equation}
\label{cycle_class}
\cl(C_u) = \Delta_u \,.
\end{equation}
\item
There~is a finite extension
field\/~$K' \supseteq K$
such that the\/
\mbox{$\bbQ$-cycle\/}
$C_u$
can be defined
over\/~$K'$.
\end{abc}
\end{prop}

\begin{defi}
The
\mbox{$\bbQ$-cycle}
$C_u$
is usually called a {\em correspondence\/} describing the action
of~$u \in E$
on~$X_{\overline{K}}$.
\end{defi}

\begin{proop}
a)
As~$[u]\colon H^2(X(\bbC), \bbQ) \to H^2(X(\bbC), \bbQ)$
is a morphism of Hodge structures, the image in
$H^4((X \!\times\! X)(\bbC), \bbC)$
under change of coefficients of the analytic Casimir~element
$\Delta^\an_u \in H^4((X \!\times\! X)(\bbC), \bbQ)$
is pure of Hodge
type~$(2,2)$.
Therefore,~the Hodge conjecture yields a
\mbox{$\bbQ$-cycle}
$\smash{C_u^\bbC}$
of
codimension~$2$~in
$$(X \!\times\! X)(\bbC) = X(\bbC) \!\times\! X(\bbC)$$
corresponding to the analytic Casimir
element~$\Delta^\an_u$.
I.e.,~of the kind~that
\begin{equation}
\label{cycle_Hodge}
\cl^\an(C_u^\bbC) = \Delta^\an_u \,.
\end{equation}
According~to the GAGA principle~\cite[Proposition~13]{Se56}, one knows that
$\smash{C_u^\bbC}$
is automatically a complex algebraic
\mbox{$\bbQ$-cycle}
on~$(X \!\times\! X)_\bbC$.
Furthermore,~the analytic cycle class map
$\cl^\an$
is compatible with the algebraic cycle class
map~$\cl$.

The~cycle class map is compatible with base change from
$\smash{\overline{K}}$
to~$\bbC$,~too.
Cf.\ \cite[tag~0FWC]{SP}. There~is, however, no need for
$\smash{C_u^\bbC}$
to descend
to~$\smash{(X \!\times\! X)_{\overline{K}}}$.
On~the other hand, algebraic equivalence is finer than homological equivalence~\cite[Proposition~19.1.1]{Fu}. Hence,~one may replace
$C_u^\bbC$
by any algebraically equivalent
\mbox{$\bbQ$-cycle}
$C_u$
on~$(X \!\times\! X)_\bbC$,
without affecting~(\ref{cycle_Hodge}).

Write
$$\smash{C_u^\bbC = r_1(C_1) + \cdots + r_m(C_m) \,,}$$
for
$\smash{r_1, \ldots, r_m \in \bbQ}$
and
$\smash{C_1, \ldots, C_m \subset (X \!\times\! X)_\bbC}$
closed subschemes of
codimension~$2$
being reduced and irreducible. We~shall treat each component individually, replacing them by algebraically equivalent ones that descend
to~$\smash{\overline{K}}$.

To~be concrete, let us first choose a projective embedding
$\smash{X \!\times\! X \hookrightarrow \Pb^N_K}$.
Then,~for each
$i$,
there is the Hilbert
scheme~$\smash{H_i := H_{P_{C_i}, X \times X}}$,
which is a projective
\mbox{$K$-scheme}
\cite[Expos\'e~221, Th\'eor\`eme~3.2]{FGA}, that parameterises closed subschemes
of~$X \!\times\! X$
having the same Hilbert polynomial
$P_{C_i}$
as~$C_i$.
The~Hilbert
scheme~$H_i$
comes equipped with a universal family
$\smash{\pi_i\colon \calC_i \hookrightarrow X \!\times\! X \!\times\! H_i \stackrel{\pr}{\twoheadrightarrow} H_i}$
that is a projective and flat morphism of~schemes. Moreover,~by \cite[Th\'eor\`eme~12.2.4.(viii)]{EGAIV}, the locus
$H_i^0 \subset H_i$,
above which the fibres
of~$\pi$
are reduced and irreducible, is~open.

The~subscheme
$\smash{C_i \subset (X \!\times\! X)_\bbC}$
gives rise to
a~\mbox{$\bbC$-rational}
point~$z_i$
on~$H_i^0$,
according to the definition of the Hilbert~scheme. In~particular, one sees that
$\smash{H_i^0 \neq \emptyset}$.
Let~$(H_i^{0})' \subseteq H_i^0$
be the irreducible component
of~$H_i^0$
containing~$z_i$.
Then~$(H_i^{0})'$
has a
\mbox{$\smash{\overline{K}}$-rational}
point~$\smash{z_i \in (H_i^{0})'(\overline{K})}$
by a weak version of Hilbert's Nullstellensatz~\cite[Corollary~13.12.i]{Ei}. The~corresponding fibre
$\smash{\widetilde{C}_i := (\calC_i \!\times_{H_i}\! \overline{\{z_i\}})_\bbC}$
is obviously algebraically equivalent
to~$C_i$,
cf.\ \cite[Example~10.3.2]{Fu}, and descends
to~$\smash{\overline{K}}$.
This~completes the proof of~a).\smallskip

\noindent
b)
immediately follows from~a).
\eop
\end{proop}

\begin{theo}[Lower bound for the field of definition of a correspondence]
\label{corr_bound}
Let\/~$X$
be a\/
$K3$~surface
with endomorphism
field\/~$E$
over a
field\/~$K\subseteq\bbC$.
Suppose~that the Hodge conjecture holds for\/
$(X \!\times\! X)(\bbC)$.
For~some\/
$u \in E$,
let\/
$C_u$
be a correspondence describing the action
of\/~$u$
on\/~$X_{\overline{K}}$
and
let\/~$K' \supseteq K$
be an extension field over which\/
$C_u$
can be~defined.\smallskip

\noindent
Then\/~$K' \supseteq E'$,
for\/~$E' := K(u)$
the subfield
of~$E$
generated
by\/~$u$.\medskip

\noindent
{\bf Proof.}
{\em
As~$C_u$
can be defined
over~$K'$,
one has
${}^\sigma (C_u) = C_u$,
for every
$\sigma \in \Gal(\overline{K}/K')$.
Since~the cycle map is
\mbox{$\smash{\Gal(\overline{K}/K')}$-equivariant},
this yields
\begin{equation}
{}^\sigma (\Delta_u) = {}^\sigma (\cl(C_u)) = \cl({}^\sigma (C_u)) = \cl(C_u) = \Delta_u \,.
\end{equation}
At~this point, Lemma~\ref{Casimir_sigma} shows in view of~(\ref{Casimir})~that
${}^{(\sigma \circ [u] \circ \sigma^{-1})} v_i = {}^{[u]} v_i$,
for
$i=1,\ldots,22$.
Thus,~the actions of
$\sigma \!\circ\! [u] \!\circ\! \sigma^{-1}$
and~$[u]$~on
$$
\bigoplus_{i=0}^2 H^i(X_{\overline{K}}, \Omega^{2-i}_{X_{\overline{K}}/\overline{K}}) \,,
$$
which are both
\mbox{$\smash{\overline{K}}$-linear}
maps,~coincide. In~other words, the action of
$\sigma$
commutes with that
of~$[u]$.

From~this, Lemma~\ref{actions_dR}.b) immediately yields that
$\sigma(u) = u$.
As~$\smash{\sigma \in \Gal(\overline{K}/K')}$
is arbitrary, this is possible only for
$u \in K'$.
I.e.,~one has
$K(u) \subseteq K'$,
as~required.%
}%
\eop
\end{theo}\pagebreak

\section{Results relying on the Hodge conjecture}
\label{Hodge_dep}

\subsubsection*{Generalities}

\begin{ttt}
Let~$X$
be a
$K3$~surface
over a field
$K\subseteq\bbC$.
Then,~as usual, a
\mbox{$\bbQ$-cycle}
$C$
of
codimension~$2$
in~$\smash{(X \!\times\! X)_{\overline{K}}}$
defines a homomorphism on algebraic de Rham cohomology~\cite[Section~3]{Kl},
$$
\gamma_C\colon \bigoplus_{i=0}^2 H^i(X_{\overline{K}}, \Omega^{2-i}_{X_{\overline{K}}/\overline{K}}) \longrightarrow \bigoplus_{i=0}^2 H^i(X_{\overline{K}}, \Omega^{2-i}_{X_{\overline{K}}/\overline{K}}) \,, \quad v \mapsto \pr_2 {}_* (\pr_1^*(v) \cup \cl(C)) \,.
$$
If~$\cl(C) = v_1^* \otimes w_1 + \cdots + v_{22}^* \otimes w_{22}$,
for
$(v_1, \ldots, v_{22})$
a basis of
$\smash{\bigoplus\limits_{i=0}^2 H^i(X_{\overline{K}}, \Omega^{2-i}_{X_{\overline{K}}/\overline{K}})}$
and
$(v_1^*, \ldots, v_{22}^*)$
the dual basis,~then
\begin{equation}
\label{motiv_map}
\gamma_C(v_i) = w_i \,,
\end{equation}
for~$i=1, \ldots, 22$.
\end{ttt}

\begin{ttt}
Let~$E$
be the endomorphism field
of~$X$.
For~each
$u \in E$,
we~extend
$[u]\colon T \to T$
to a
\mbox{$\bbQ$-linear}
map
$$[u]\colon H^2(X(\bbC), \bbQ) \to H^2(X(\bbC), \bbQ)$$
by mapping
$H_\alg$
identically
to~$0$.
This~is a particular case of the extensions considered in the section~above.

Suppose,~in addition, that the Hodge conjecture holds
for~$(X \!\times\! X)(\bbC)$.
Then,~according to Proposition~\ref{corr}.a), there is a
\mbox{$\bbQ$-cycle}
$C_u$
of
codimension~$2$
in~$\smash{(X \!\times\! X)_{\overline{K}}}$
satisfying the condition
that~$\cl(C_u) = \Delta_u$.
By~(\ref{motiv_map}), for~every
$\smash{u \in \bigoplus_{i=0}^2 H^i(X_{\overline{K}}, \Omega^{2-i}_{X_{\overline{K}}/\overline{K}})}$,
one~has
$$\gamma_{C_u}(v) = {}^{[u]} v \,.$$

Furthermore,~for arbitrary
$\sigma\in\Gal(\overline{K}/K)$,
let us consider the
\mbox{$\bbQ$-cycle}
${}^\sigma (C_u)$,
which is again of
codimension~$2$
in~$\smash{(X \!\times\! X)_{\overline{K}}}$.
One~clearly has
$\cl({}^\sigma (C_u)) = {}^\sigma (\Delta_u)$.
Thus,~from Lemma~\ref{Casimir_sigma}, together with~(\ref{motiv_map}), one concludes that
\begin{equation}
\label{motiv_map_konj}
\gamma_{{}^\sigma (C_u)}(v) = {}^{\sigma \circ [u] \circ \sigma^{-1}} v \,,
\end{equation}
for every
$\smash{u \in \bigoplus\limits_{i=0}^2 H^i(X_{\overline{K}}, \Omega^{2-i}_{X_{\overline{K}}/\overline{K}})}$.
\end{ttt}

\begin{ttt}
Analogously~to the above, a
\mbox{$\bbQ$-cycle}
$C$
of
codimension~$2$
on~$\smash{(X \!\times\! X)_{\overline{K}}}$
also defines a homomorphism on
\mbox{$l$-adic}
cohomology,
$$
\gamma_{C,l}\colon H^2_\et(X_{\overline{K}}, \overline\bbQ_l(1)) \longrightarrow H^2_\et(X_{\overline{K}}, \overline\bbQ_l(1)) \,, \quad v \mapsto \pr_2 {}_* (\pr_1^*(v) \cup \cl_l(C)) \,.
$$

Suppose~that
$k = K$
is a number~field. Then~there is the comparison
isomorphism~$\iota_{X,\frakl}$,
which is compatible with the cycle map, cup products, and the K\"unneth decomposition \cite[Theorem~II.3.1, cf.\ Paragraph \ref{nf_case}]{Fa}. So one has a commutative~diagram
$$
\xymatrix{
H^2_\et(X_{\overline{k}}, \bbC_l(1)) \ar@{->}[d]_\cong^{\iota_{X,\frakl}} \ar@{->}[r]^{\gamma_{C,l} \otimes \id} & H^2_\et(X_{\overline{k}}, \bbC_l(1)) \ar@{->}[d]_\cong^{\iota_{X,\frakl}} \\
\smash{\bigoplus\limits_{i=0}^2}\, H^i(X_{\overline{k}},\Omega^{2-i}_{X_{\overline{k}}/\overline{k}}) \!\otimes_{\overline{k}}\! \bbC_l(i\!-\!1)\, \ar@{->}[r]^{\gamma_C \otimes \id} & \,\smash{\bigoplus\limits_{i=0}^2}\, H^i(X_{\overline{k}},\Omega^{2-i}_{X_{\overline{k}}/\overline{k}}) \!\otimes_{\overline{k}}\! \bbC_l(i\!-\!1) \,.
}
$$
Therefore,~Proposition~\ref{transf} shows~that
\begin{equation}
\label{cyc_ladic}
\gamma_{C_u,l}(v) = {}^{[u]} v \,,
\end{equation}
for every
$v \in T_{\overline{l}}$,~too.
Moreover,~the
\mbox{$l$-adic}
cycle map
$\cl_l$,
as well as cup products and the K\"unneth decomposition, are compatible with the action
of~$\Gal(\overline{k}/k)$.
Thus,~(\ref{cyc_ladic}) implies that
$\gamma_{{}^\sigma (C_u),l}({}^\sigma v) = {}^\sigma ({}^{[u]} v)$,
which is equivalent~to
\begin{equation}
\label{cyc_ladic_sigma}
\gamma_{{}^\sigma (C_u),l}(v) = {}^{\sigma \circ [u] \circ \sigma^{-1}} v \,.
\end{equation}
\end{ttt}\vskip2mm

\begin{lem}
\label{Hodge_corr}
Let\/~$X$
be a\/
$K3$~surface
with endomorphism
field\/~$E$
over a
field\/
$K\subseteq\bbC$
and\/
$C$
a\/
\mbox{$\bbQ$-cycle}
of
codimension\/~$2$
in\/~$\smash{(X \!\times\! X)_{\overline{K}}}$.
Furthermore,~let some nonzero\/
$\smash{\omega \in H^0(X_{\overline{K}}, \Omega^2_{X_{\overline{K}}/\overline{K}})}$
be~given.

\begin{abc}
\item
Then\/~$\gamma_C(\omega) = u \cdot \omega$,
for some scalar
factor\/~$u$
that is necessarily an element
of\/~$E$.
\item
Suppose,~in particular, that\/
$\gamma_C(\omega) = 0$.
Then~the homomorphism\/
$\gamma_C$
is identically zero on\/
$\smash{H^i(X_{\overline{K}}, \Omega^{2-i}_{X_{\overline{K}}/\overline{K}})}$,
for\/
$i = 0,2$,
as well as
on\/~$\smash{V_{\overline{K}}}$.
\end{abc}\smallskip

\noindent
{\bf Proof.}
{\em
The~homomorphism
$\gamma_C^\an\colon H^2(X(\bbC), \bbQ) \to H^2(X(\bbC), \bbQ)$
on complex cohomology induced
by~$\gamma_C$
may be described by
$\gamma_C^\an\colon v \mapsto \pr_2 {}_* (\pr_1^*(v) \cup \cl^\an(C_\bbC))$,
for~$C_\bbC$
the complexification
of~$C$.
Therefore,~$\gamma_C^\an$
is a morphism of Hodge structures~\cite[Section~2.5]{GH}, which implies that it must map the transcendental part
$T \subset H^2(X(\bbC), \bbQ)$
to~itself.
Clearly,~$\gamma_C^\an |_T \colon T \to T$
is an endomorphism of the Hodge
structure~$T$.
I.e.,~$\gamma_C^\an |_T$
and~hence
$$\textstyle \gamma_C |_{\!\bigoplus\limits_{i=0,2}\! H^i(X_{\overline{K}}, \Omega^{2-i}_{X_{\overline{K}}/\overline{K}}) \oplus V_{\overline{K}}} \colon \bigoplus\limits_{i=0,2}\! H^i(X_{\overline{K}}, \Omega^{2-i}_{X_{\overline{K}}/\overline{K}}) \oplus V_{\overline{K}} \longrightarrow \bigoplus\limits_{i=0,2}\! H^i(X_{\overline{K}}, \Omega^{2-i}_{X_{\overline{K}}/\overline{K}}) \oplus V_{\overline{K}} $$
itself, too, must be given by an element
$u \in E$.\smallskip

\noindent
a)
In~particular, one has that
$\smash{\gamma_C |_{H^0(X_{\overline{K}},\Omega^2_{X_{\overline{K}}/\overline{K}})}}$
is simply the map
$\gamma_C\colon \omega \mapsto u \cdot \omega$,
which proves the~assertion.\smallskip

\noindent
b)
Here,~the assumption yields that
$u=0$.
Therefore,
$\smash{\gamma_C |_{\!\bigoplus\limits_{i=0,2}\! H^i(X_{\overline{K}}, \Omega^{2-i}_{X_{\overline{K}}/\overline{K}}) \oplus V_{\overline{K}}} = 0}$,
as~claimed.
}
\eop
\end{lem}

\subsubsection*{Normality of the endomorphism field}
\leavevmode\medskip

\noindent
The result below is well-known in the context of CM elliptic curves~\cite[Proposition~5.17.(3)]{Sh}.

\begin{theo}[Normality of the endomorphism field]
\label{normality}
Let\/~$X$
be a\/
$K3$~surface
with endomorphism
field\/~$E$
over a
field\/
$K\subseteq\bbC$.
Suppose~that the Hodge conjecture holds
for\/~$(X \!\times\! X)(\bbC)$.
Then~the composite field\/
$KE$
is normal
over\/~$K$.\medskip

\noindent
{\bf Proof.}
{\em
Let~$u_0 \in E$
be a primitive~element. As~the Hodge conjecture is assumed to hold
for~$(X \!\times\! X)(\bbC)$,
by Proposition~\ref{corr}.a), there is a
\mbox{$\bbQ$-cycle}
$C_{u_0}$
of
codimension~$2$
in~$\smash{(X \!\times\! X)_{\overline{K}}}$
satisfying
$\cl(C_{u_0}) = \Delta_{u_0}$.
For~an arbitrary
$\smash{\sigma \in \Gal(\overline{K}/K)}$,
let us consider the
\mbox{$\bbQ$-cycle}
${}^\sigma (C_{u_0})$,
which is again of
codimension~$2$
in~$\smash{(X \!\times\! X)_{\overline{K}}}$.
Concerning~the homomorphism on algebraic de Rham cohomology defined
by~${}^\sigma (C_{u_0})$,
one knows from formula~(\ref{motiv_map_konj})~that
$$\gamma_{{}^\sigma (C_{u_0})}(v) =  {}^{\sigma \circ [u_0] \circ \sigma^{-1}} v \,,$$
for every
$\smash{u \in \bigoplus_{i=0}^2 H^i(X_{\overline{K}}, \Omega^{2-i}_{X_{\overline{K}}/\overline{K}})}$.

For~$\smash{v \in H^0(X_{\overline{K}},\Omega^2_{X_{\overline{K}}/\overline{K}})}$,
Lemma~\ref{actions_dR}.a) computes the right hand side~explicitly. One~has
\begin{equation}
\label{v1}
\gamma_{{}^\sigma (C_{u_0})}(v) = \sigma^{-1}(u_0) \cdot v \,.
\end{equation}
At~this point, from Lemma~\ref{Hodge_corr}.a), one sees that
$\sigma^{-1}(u_0) \in E \subseteq KE$
is~enforced.
Since~$u_0 \in KE$
is a primitive element relative
to~$K$
and
$\smash{\sigma \in \Gal(\overline{K}/K)}$
is arbitrary, this shows that
$KE$
is normal
over~$K$,
which completes the~proof.
}
\eop
\end{theo}

\subsubsection*{The splitting field and the jump character in the CM case}

\begin{lem}[Commuting actions]
\label{comm-dR}
Let\/~$X$
be a\/
$K3$~surface
with endomorphism
field\/~$E$
over a
field\/
$K\subseteq\bbC$.
Suppose~that the Hodge conjecture holds
for\/~$(X \!\times\! X)(\bbC)$.

\begin{abc}
\item
Then,~for every\/
$u \in E$
and\/
$\sigma \in \Gal(\overline{K}/K) \subseteq \Aut_\bbQ(\overline{K})$,
one has\/
$${}^{\sigma \circ [u] \circ \sigma^{-1}} \eta = {}^{[\sigma^{-1}(u)]} \eta \,,$$
for arbitrary\/
$\smash{\eta \in \bigoplus_{i=0}^2 H^i(X_{\overline{K}}, \Omega^j_{X_{\overline{K}}/\overline{K}})}$.
\item
Suppose,~in particular, that\/
$k = K$
is a number~field.
Let\/~$l$
be any prime number and
let\/~$v \in T_{\overline{l}}$
be~arbitrary. Then,~for every\/
$\smash{\sigma \in \Gal(\overline{k}/k)}$
and
each\/~$u \in E$,
$${}^{\sigma \circ [u] \circ \sigma^{-1}} v = {}^{[\sigma^{-1}(u)]} v \,.$$
\end{abc}

\noindent
{\bf Proof.}
{\em
a)
For~the
\mbox{$\bbQ$-cycle}
$D := {}^\sigma (C_u) - C_{[\sigma^{-1}(u)]}$
on~$(X \!\times\! X)_{\overline{K}}$,
one has
$$\gamma_D(\eta) = {}^{\sigma \circ [u] \circ \sigma^{-1}} \eta - {}^{[\sigma^{-1}(u)]} \eta \,,$$
for every
$\smash{\eta \in \bigoplus_{i=0}^2 H^i(X_{\overline{K}}, \Omega^{2-i}_{X_{\overline{K}}/\overline{K}})}$.
Therefore,~Lemma~\ref{actions_dR} shows that
$\gamma_D$
vanishes on
$\smash{H^0(X_{\overline{K}}, \Omega^2_{X_{\overline{K}}/\overline{K}})}$.
Hence,~according to Lemma~\ref{Hodge_corr}.b),
$\gamma_D$
is the zero map on the whole~of
$$\textstyle \bigoplus\limits_{i=0,2}\! H^i(X_{\overline{K}}, \Omega^{2-i}_{X_{\overline{K}}/\overline{K}}) \oplus V_{\overline{K}} \,.$$
As~the actions of
$u$
and~$\sigma^{-1}(u)$
are assumed to be identically zero on the algebraic part of the cohomology,
$\gamma_D$
is the zero map altogether, which completes the proof of~a).\smallskip

\noindent
b)
In~view of~(\ref{motiv_map}), the result of~a) shows that the
\mbox{$\bbQ$-cycle}
$D$
on~$(X \!\times\! X)_{\overline{k}}$
is homologically equivalent to~zero. This~property holds for
\mbox{$l$-adic}
cohomology, too, which implies the assertion, due to formulae (\ref{cyc_ladic}) and~(\ref{cyc_ladic_sigma}). Note~here that clearly
$D$~is
numerically equivalent to~zero.
As~$D$
is
a cycle of
codimension~$2$,
this is known to imply homological equivalence to zero with respect to any Weil cohomology theory~\cite[Corollary~1]{Lie}.%
}%
\eop
\end{lem}

\begin{theo}[The splitting field and the jump character]
\label{jump_CM_alg}
Let\/~$X$
be a\/
$K3$~surface
over a number
field\/~$k$.
Assume~that the endomorphism
field\/~$E$
of\/~$X$
is a CM~field.

\begin{abc}
\item
Then,~for every prime
number\/~$l$,
there is a natural isomorphism
$$G_{X,\overline{l}}/G_{X,\overline{l}}^0 \stackrel{\cong}{\longrightarrow} \Gal(kE/k) \,.$$
I.e.,~the splitting field
of\/~$\smash{G_{X,\overline{l}}/G_{X,\overline{l}}^0}$
is
exactly\/~$kE$.
\item
The~jump character\/
$\tau_X\colon \Gal(\overline{k}/k) \to \{\pm1\}$
is
$$
\sigma \mapsto \left\{
\begin{array}{rrr}
1 \,,                                                                             & \hspace{.7cm}\textstyle\text{~if\/~} \smash{\frac{22-\rk\Pic X_{\overline{k}}}{[E:\bbQ]}} \text{~is~even} \,, \\
\text{sign~of~the~natural~permutation~}                                           & \\[-.4mm]
\text{action~of~} \sigma \in \smash{\Gal(\overline{k}/k)} \text{~on~the~}         & \\[-.4mm]
\text{conjugates\/~} u_1, \overline{u}_1, \ldots, u_d, \overline{u}_d \text{~of~} & \\[-.4mm]
\text{a~primitive~element~} u_0 \in E \,,                                         & \textstyle\text{~if\/~} \smash{\frac{22-\rk\Pic X_{\overline{k}}}{[E:\bbQ]}} \text{~is~odd} \,. \;\;
\end{array}
\right.
$$
\end{abc}\smallskip

\noindent
{\bf Proof.}
{\em
It~is known that the Hodge conjecture holds for
$(X \!\times\! X)(\bbC)$
\cite[Theorem 5.4]{RM}. Thus, in view of Lemma~\ref{comm-dR}.b), both assertions follow directly from the fundamental Lemma~\ref{fund_lem}.
}
\eop
\end{theo}

\begin{coro}
\label{jump}
Let\/
$X$
be a\/
$K3$~surface
over a number
field\/~$k$.
Suppose~that the endomorphism field
of\/~$X$
is an imaginary quadratic field
$E = \bbQ(\sqrt{-\delta})$,
for\/~$\delta\in\bbN$.\smallskip

\noindent
Then~the jump character\/
$\tau_X\colon \Gal(\overline{k}/k) \to \{\pm1\}$
is given~by
$$
\Frob_\frakp \mapsto \left\{
\begin{array}{cr}
 1 \,, \;\;\;                           & {\it ~if} \rk\Pic X_{\overline{k}} \equiv 2 \pmodulo 4 \,, \\[-.3mm]
\textstyle (\frac{-\delta}\frakp) \,, \;\;\; & {\it ~if} \rk\Pic X_{\overline{k}} \equiv 0 \pmodulo 4 \,.
\end{array}
\right.
$$
\end{coro}

\section{Examples}

\subsubsection*{Complex multiplication}

\begin{ex}
Let~$X'_1$
be the double cover
of~$\Pb^2_\bbQ$,
given by
$$w^2 = xyz (x+y+z) (x+2y+3z) (5x+8y+20z)$$
and
$X_1$
the
$K3$~surface
obtained as the minimal desingularisation
of~$X_1'$.
Then, as shown in \cite[Example~5.7]{EJ22b}, the geometric Picard rank
of~$X_1$
is~$16$
and the endomorphism field
of~$X_1$
is~$E = \bbQ(i)$.
Thus,~Theorem~\ref{jump_CM_alg}.a) implies that
$\smash{G_{X_1,\overline{l}}/G_{X_1,\overline{l}}^0 \cong \Gal(E/\bbQ) = \bbZ/2\bbZ}$
and Corollary~\ref{jump} shows that the jump character
$\tau_{X_1}$
is given
by~$(\frac{-1}.)$.
Both facts have been obtained before, using other~methods.
\end{ex}

\begin{ex}
\label{Kummer_cyc}
Let~$X_2 := \Kum(\J(C))$
be the Kummer surface associated with the Jacobian of the
genus~$2$
curve~$C$
over~$\bbQ$,
given by
$w^2 = x^5-1$.

\begin{abc}
\item
Then~the geometric Picard rank
of~$X_2$
is~$18$.
\item
The endomorphism field
of~$X_2$
is~$E = \bbQ(\zeta_5)$
and the jump character
$\tau_{X_2}$
is given
by~$(\frac5.)$.
\end{abc}\smallskip

\noindent
{\bf Proof.}
a)
According~to~\cite[\S7]{vW}, the Jacobian
$\J(C)$
of~$C$
has CM by the quartic field
$\bbQ(\zeta_5)$.
This~yields that
$\smash{\rk\NS\J(C)_\bbC = [\bbQ(\zeta_5)\!:\!\bbQ]/2 = 2}$~\cite[Section 21, Application III]{Mu}. Finally,~it is well-known (cf.\ \cite{Lim} or \cite[Fact~4.1]{EJ12})~that
$$\smash{\rk\Pic X_{2,\bbC} = \rk\Pic\Kum(\J(C))_\bbC = \rk\NS\J(C)_\bbC + 16 \,.}$$

\noindent
b)
For~the transcendental part
$T \subset H^2(X_2(\bbC), \bbQ)$
of the cohomology, one has an inclusion
$\smash{T \hookrightarrow H^2(\J(C)(\bbC), \bbQ) = \Exterior^2 H^1(C(\bbC), \bbQ)}$.
On~the right hand side, let
$\zeta_5$
act as
$\smash{v \wedge v' \mapsto {}^{[\zeta_5]} v \wedge {}^{[\zeta_5]} v'}$.
This~is well-defined and a morphism of Hodge~structures.
Moreover,
on~$T$,
the minimal polynomial of this morphism is
$\smash{(T^5-1)/(T-1)}$,
so there is an extension to an action of the whole
of~$\bbQ(\zeta_5)$.
I.e.,~one has
$E \supseteq \bbQ(\zeta_5)$
and the inclusion the other way round follows from the fact that
$T$
carries the structure of an
\mbox{$E$-vector}
space, which implies
$\smash{[E\!:\!\bbQ] \mid (22-\rk\Pic X_{2,\overline\bbQ})}$.
The~final claim follows from Theorem~\ref{jump_CM_alg}.b). It~may as well be obtained using \cite[Algorithm~2.6.1]{CEJ}.\eop\medskip

In~view of these results, Proposition~\ref{CM_general}.a) shows that
$\smash{G_{X_2,\overline{l}}^0 \cong [\overline\bbQ_l^*]^2}$.
Moreover,
$\smash{G_{X_2,\overline{l}}/G_{X_2,\overline{l}}^0 \cong \Gal(E/\bbQ) = \bbZ/4\bbZ}$,\vspace{.15mm}
by Theorem~\ref{jump_CM_alg}.a). Based~on point counting
on~$C$,
which is faster than counting points
on~$X_2$,
we determined
$\smash{\Tr(\varrho_{X_2,\overline{l}}(\Frob_p))}$,
for all prime numbers
$p \leq 5 \!\cdot\! 10^8$.
If~$p \not\equiv 1 \pmod 5$
then
$\smash{\Tr(\varrho_{X_2,\overline{l}}(\Frob_p)) = 0}$.
The~experimental distribution of the traces is shown in the histogram below, plotted against the theoretical distribution, according to the Sato--Tate conjecture, the density of which is given
by~$\K(1-t^2/16)/8\pi^2$,
for~$\K$
the complete elliptic integral of the first~kind. Cf.\ \cite[Section~2 and formula~(6)]{EJ22b}. The~histogram ignores about the primes
$p \not\equiv 1 \pmod 5$,
which would add a spike of
mass~$3/4$
above~$0$.\medskip

\begin{figure}[H]
\begin{center}
\includegraphics[scale=0.4]{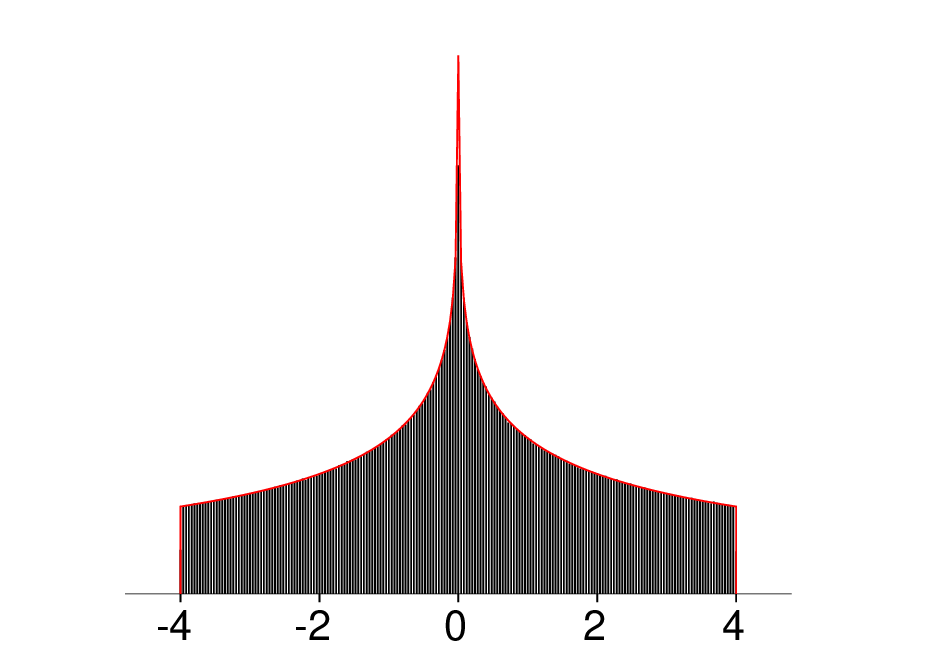}
\end{center}\vskip-5mm
\caption{Theoretical and experimental trace distributions for the
$K3$~surface
$X_2$
over
$k=\bbQ$
of geometric Picard
rank~$18$
having CM by
$\smash{E=\bbQ(\zeta_5)}$
}
\end{figure}
\end{ex}

\begin{ex}
Let~$X_3$
be the Kummer surface associated with the split abelian surface
$E_1 \!\times\! E_2$,
for
$E_1$
and~$E_2$
the elliptic curves 
over~$\bbQ$,
given by
$E_1\colon y^2 = x^3+x$
and
$E_2\colon y^2 = x^3+4x^2+2x$.

\begin{abc}
\item
Then~the geometric Picard rank
of~$X_3$
is~$18$.
\item
The endomorphism field
of~$X_3$
is~$E = \bbQ(\sqrt{2},i)$
and the jump character
$\tau_{X_3}$
is~trivial.
\end{abc}\smallskip

\noindent
{\bf Proof.}
a)
One has
$j(E_1) = 12^3$
and~$j(E_2) = 20^3$,
hence the curves have CM by
$\bbQ(i)$
and~$\smash{\bbQ(\sqrt{-2})}$,
respectively~\cite[\S12.A]{Co}. This~yields that
$\smash{\End((E_1 \!\times\! E_2)_\bbC) = \bbQ(\sqrt{2},i)}$
and that
$\smash{\rk\NS((E_1 \!\times\! E_2)_\bbC) = [\bbQ(\sqrt{2},i)\!:\!\bbQ]/2 = 2}$~\cite[Section 21, Application III]{Mu},
which implies the claim.\smallskip

\noindent
b)
For~$T \subset H^2(X_3(\bbC), \bbQ)$,
the K\"unneth formula yields a natural isomorphism
$T \cong H^1(E_1(\bbC), \bbQ) \!\otimes_\bbQ\! H^1(E_2(\bbC), \bbQ)$.
Moreover, the action
of~$\smash{\bbQ(\sqrt{2},i)}$
on~$T$
is given by
$\smash{{}^{[i]} (v \!\otimes\! v') := {}^{[i]} v \!\otimes\! v'}$
and
$\smash{{}^{[\sqrt{-2}]} (v \!\otimes\! v') := v \!\otimes\! {}^{[\sqrt{-2}]} v'}$.
Again,~the equality
$\smash{E = \bbQ(\sqrt{2},i)}$
is enforced by the fact that
$T$
carries the structure of an
\mbox{$E$-vector}
space. The~final claim is obtained using \cite[Algorithm~2.6.1]{CEJ}.\eop\medskip

Note~that the result on the jump character as well follows from Theorem~\ref{jump_CM_alg}.b). Based~on point counting on
$E_1$
and~$E_2$,
we determined
$\smash{\Tr(\varrho_{X_3,\overline{l}}(\Frob_p))}$,
for all prime numbers
$\smash{p \leq 5 \!\cdot\! 10^8}$.
If~$p \not\equiv 1 \pmod 8$
then
$\smash{\Tr(\varrho_{X_3,\overline{l}}(\Frob_p)) = 0}$.
The~experimental distribution of the traces is shown in the histogram above, plotted against the theoretical distribution, according to the Sato--Tate conjecture, which is the same as in Example~\ref{Kummer_cyc}. Note~here that one~has\vspace{.3mm}
$\smash{G_{X_3,\overline{l}}^0 \cong [\overline\bbQ_l^*]^2}$,
again. Furthermore, Theorem~\ref{jump_CM_alg}.a) shows that
$\smash{G_{X_3,\overline{l}}/G_{X_3,\overline{l}}^0 \cong \Gal(E/\bbQ) = [\bbZ/2\bbZ]^2}$.
Again,~the spike
above~$0$
of
mass~$3/4$
is~ignored.\medskip

\begin{figure}[H]
\begin{center}
\includegraphics[scale=0.4]{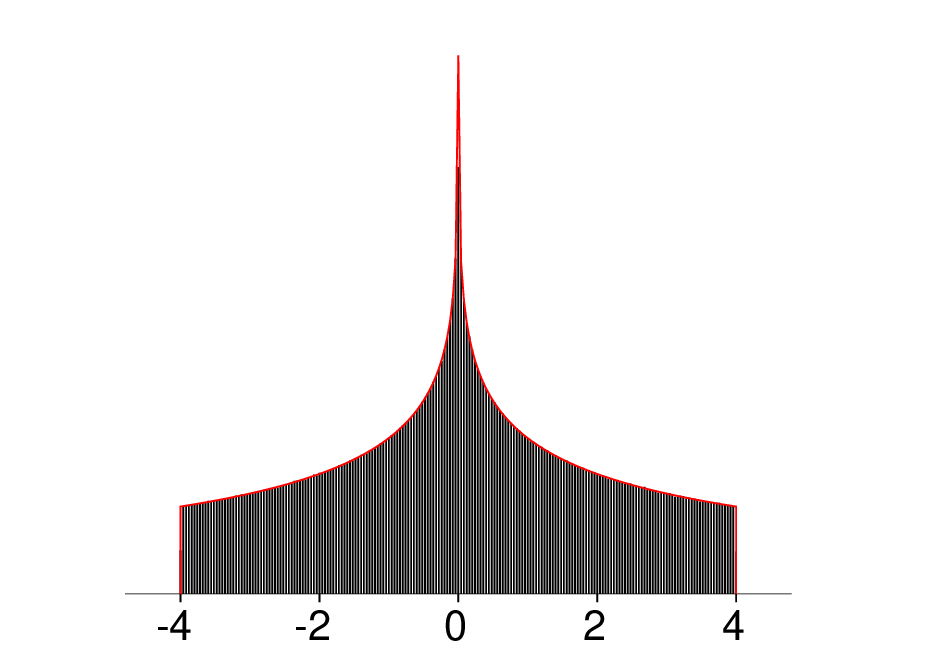}
\end{center}\vskip-5mm
\caption{Theoretical and experimental trace distributions for the
$K3$~surface
$X_3$
over
$k=\bbQ$
of geometric Picard
rank~$18$
having CM by
$\smash{E=\bbQ(\sqrt{2},i)}$
}
\end{figure}
\end{ex}

\begin{ex}
Let~$X_4'$
be the double cover
of~$\Pb^2_\bbQ$,
given by
$$w^2 = xyz (x^3 - 3x^2z - 3xy^2 - 3xyz + y^3 + 9y^2z + 6yz^2 + z^3)$$
and
$X_4$
the
$K3$~surface
obtained as the minimal desingularisation
of~$X_4'$.
Then, as shown in \cite[Example~5.11]{EJ22b}, the geometric Picard rank
of~$X_4$
is~$16$.
Moreover,~there is strong evidence that the endomorphism field
is~$E = \bbQ(\zeta_9 + \zeta_9^{-1},\sqrt{-1})$.
Assuming~this, Theorem~\ref{jump_CM_alg}.a) implies~that
$$\smash{G_{X_4,\overline{l}}/G_{X_4,\overline{l}}^0 \cong \Gal(E/\bbQ) = \bbZ/6\bbZ}$$
and Theorem~\ref{jump_CM_alg}.b) shows that the jump character
$\tau_{X_4}$
is given
by~$\smash{(\frac{-1}.)}$.
The~latter fact has been obtained~before. For~the former, strong evidence has been~reported.
\end{ex}

\begin{ex}
Put
$k := \bbQ[t]/(t^3 - t^2 - 4t + 1)$
and let
$X_5'$
be the double cover
of~$\Pb^2_k$,
given by
$$w^2 = xyz(x+y+z)(x + \alpha y + \beta z)(x + \gamma y + \delta z) \,,$$
for
$\alpha := \overline{\frac{-26t^2 - 23t + 16}9}$,
$\beta := \overline{\frac{-61t^2 + 125t + 95}{121}}$,
$\gamma := \overline{\frac{-t^2 - 4t + 11}9}$,
and
$\delta := \overline{\frac{-46t^2 + 5t + 149}{121}}$.
Let~$X_5$
be the
$K3$~surface
obtained as the minimal desingularisation
of~$X_5'$.
Then~the geometric Picard rank
of~$X_5$
is~$16$.\medskip

\noindent
{\bf Proof.}
An upper bound of 16 is provided by the reduction
$X_{5,\frakp}$
at the unique prime
ideal~$\frakp$
of~$k$
of ideal
norm~$29$,
which is
of geometric Picard
rank~$16$.
On~the other hand, the ramification locus has 15 singular~points. Thus,
$X_5$
contains 15
\mbox{$(-2)$-curves},
which are mutually skew, and hence provide a lower bound
of~$16$,
together with the inverse image of a general line
on~$\Pb^2_k$.\eop\medskip

There~is strong evidence that the endomorphism field
of~$X_5$
is
$E = k(i)$.
Note~that the cubic field
$k$~is
totally real, so
$E$~is
indeed a CM~field. Moreover,
$E$
is not normal
over~$\bbQ$,
but, in agreement with Theorem~\ref{normality}, one has that
$kE = k(i)$
is normal
over~$k$.
If~this is true then Proposition~\ref{CM_general}.a) shows that
$\smash{G_{X_5,\overline{l}}^0 \cong [\overline\bbQ_l^*]^3}$.
Furthermore,
$\smash{G_{X_5,\overline{l}}/G_{X_5,\overline{l}}^0 \cong \Gal(kE/k) = \Gal(k(i)/k) \cong \bbZ/2\bbZ}$
and the jump character
$\tau_{X_5}$
is given
\!by~$\smash{\!(\frac{-1}.)}$.

The evidence is as~follows. We~calculated the characteristic polynomial
of~$\Frob_\frakp$
on~$\smash{T_{\overline{l}} \subset H^2_\et(X_{5,\overline{k}}, \overline\bbQ_l(1))}$,
for all prime
ideals~$\frakp$
of~$k$
of ideal
norm~$< \!10\,000$,
at which
$X_5$
has good~reduction. It~turned out that, for
each~$\frakp$,
$\rk\Pic X_{5,\overline\frakp}$
is either
$16$
or~$22$.
More~precisely, if
$\frakp$
is inert
in~$k(i)$
then one always has that
$\smash{\Tr(\varrho_{X_5,\overline{l}}(\Frob_\frakp)) = 0}$
and that the reduction is of geometric Picard
rank~$22$.
On~the other hand, if
$\frakp$
splits
in~$k(i)$
then the geometric Picard rank of the reduction happens to
be~$16$,
each~time. Moreover, the characteristic polynomial
of~$\Frob_\frakp$
on~$T_{\overline{l}}$
splits off two linear factors
over~$k(i)$.
Over~$\smash{k^{(n)}(i)}$,
it splits into linear factors completely. 
If~$\frakp$
is a prime of
degree~$3$
then characteristic polynomial
of~$\Frob_\frakp$
on~$T_{\overline{l}}$
is a perfect cube of a quadratic polynomial. Furthermore, the height~\cite{AM} of the reduction
$X_{5,\frakp}$
coincides with the degree
of~$\frakp$.
In~particular,
$X_{5,\frakp}$
is ordinary if and only if
$\frakp$~is
a degree one prime ideal that splits
in~$k(i)$.

In~addition, we calculated the trace
of~$\smash{\Frob_\frakp}$,
for all
primes~$\frakp$
up to ideal
norm~$10^7$.
The observation that
$\smash{\Tr(\varrho_{X_5,\overline{l}}(\Frob_\frakp)) = 0}$,
for every inert
prime~$\frakp$,
extends up to this bound. Moreover, for each split prime, one has exactly
$$\nu_p(\Tr(\varrho_{X_5,\overline{l}}(\Frob_\frakp))) = [\calO_k/\frakp : \bbF_{\!p}] - 1 \,,$$
for
$\nu_p\colon \bbQ^* \to \bbZ$
the normalised
\mbox{$p$-adic}
valuation and
$p$
the prime number lying
below~$\frakp$.
The information is summarised in the table~below. Note,~in particular, that the observation concerning ordinarity made above extends up to the bound
of~$10^7$.

\begin{table}[H]
\begin{center}
{\footnotesize
\begin{tabular}{|l||c|c|}
\hline
\diagbox[width=3.8cm]{Degree of $\frakp$}{Behaviour of $\frakp$ \\ in $k(i)/k$} & inert & split \\\hline\hline
$1$ & $\Tr(\varrho_{X_5,\overline{l}}(\Frob_\frakp)) = 0$ & $\nu_p(\Tr(\varrho_{X_5,\overline{l}}(\Frob_\frakp))) = 0$ \\\hline
$2$ & not possible                             & $\nu_p(\Tr(\varrho_{X_5,\overline{l}}(\Frob_\frakp))) = 1$ \\\hline
$3$ & $\Tr(\varrho_{X_5,\overline{l}}(\Frob_\frakp)) = 0$ & $\nu_p(\Tr(\varrho_{X_5,\overline{l}}(\Frob_\frakp))) = 2$ \\\hline
\end{tabular}
}
\end{center}\vspace{-2mm}

\caption{The traces of
$\Frob_\frakp$
for the
$K3$~surface
$X_5$
over a non-normal cubic number field
$k$
of geometric Picard
rank~$16$
having CM by
$\smash{E=k(i)}$}\vspace{-2mm}
\end{table}

\begin{rem}
Having~chosen a suitable basis
of~$H^2(X_5(\bbC), \bbZ)$,
as described in \cite[Corollary~3.14]{EJ22a}, the restricted period point \cite[Definition~3.15]{EJ22a}
of
$X_5$~is
\begin{equation}
\label{periods}
(1 : \overline{t^2\!-\!t\!-\!2} : \overline{1\!-\!t} : i : \overline{t^2\!-\!t\!-\!2}\!\cdot\!i : \overline{1\!-\!t}\!\cdot\!i) \,,
\end{equation}
up to an error of less
than~$10^{-229}$.
Here,
$k$~is
considered as a subfield
of~$\bbR$,
via the embedding, given
by~$\overline{t} \mapsto 0.23912\ldots\;$
For~the
\mbox{$\bbQ$-linear}
action
of~$k(i)$,
given~by
\[
\overline{t} \mapsto
\left(\!\!
\begin{smallmatrix}
 \phantom{-}1 & \phantom{-}1 & -2 & \phantom{-}0 & \phantom{-}0 & \phantom{-}0 \\
 \phantom{-}0 & \phantom{-}0 & -1 & \phantom{-}0 & \phantom{-}0 & \phantom{-}0 \\
-1 & -2 & \phantom{-}0 & \phantom{-}0 & \phantom{-}0 & \phantom{-}0 \\
\phantom{-}0 & \phantom{-}0 & \phantom{-}0 & \phantom{-}1 & \phantom{-}1 & -2 \\
\phantom{-}0 & \phantom{-}0 & \phantom{-}0 & \phantom{-}0 & \phantom{-}0 & -1 \\
\phantom{-}0 & \phantom{-}0 & \phantom{-}0 & -1 & -2 & \phantom{-}0
\end{smallmatrix}
\right)
\quad\;\text{and}\qquad
i \mapsto
\left(\!\!
\begin{smallmatrix}
 \phantom{-}0 & \phantom{-}0 & \phantom{-}0 & \phantom{-}1 & \phantom{-}0 & \phantom{-}0 \\
 \phantom{-}0 & \phantom{-}0 & \phantom{-}0 & \phantom{-}0 & \phantom{-}1 & \phantom{-}0 \\
 \phantom{-}0 & \phantom{-}0 & \phantom{-}0 & \phantom{-}0 & \phantom{-}0 & \phantom{-}1 \\
 -1 & \phantom{-}0 & \phantom{-}0 & \phantom{-}0 & \phantom{-}0 & \phantom{-}0 \\
 \phantom{-}0 & -1 & \phantom{-}0 & \phantom{-}0 & \phantom{-}0 & \phantom{-}0 \\
 \phantom{-}0 & \phantom{-}0 & -1 & \phantom{-}0 & \phantom{-}0 & \phantom{-}0
\end{smallmatrix}
\right)
,
\]
the one-dimensional vector space underlying~(\ref{periods}) is an eigenspace, which would prove CM
by~$k(i)$
if we knew the period point exactly \cite[Theorem~6.29.b)]{EJ22a}. The~search for this example has actually been based on period integration \cite[Algorithm~5.2]{EJ22a}.
\end{rem}

The~experimental distribution of the traces for all prime ideal
$\frakp$
in~$k$
of ideal norm
$< \!10^7$
is shown in the histogram below, plotted against the theoretical distribution, according to the Sato--Tate conjecture. The spike of
mass~$1/2$
above~$0$
is not~shown. The~density of the theoretical distribution agrees, up to scaling, with the one discussed in~\cite[Section~3, second of Examples~B]{EJ22b}. There~is an explicit formula, which is, however, rather complicated \cite[formula~(7)]{EJ22b}.\vskip-6.5mm

\begin{figure}[H]
\label{CM_rel_plot}
\begin{center}
\includegraphics[scale=0.4]{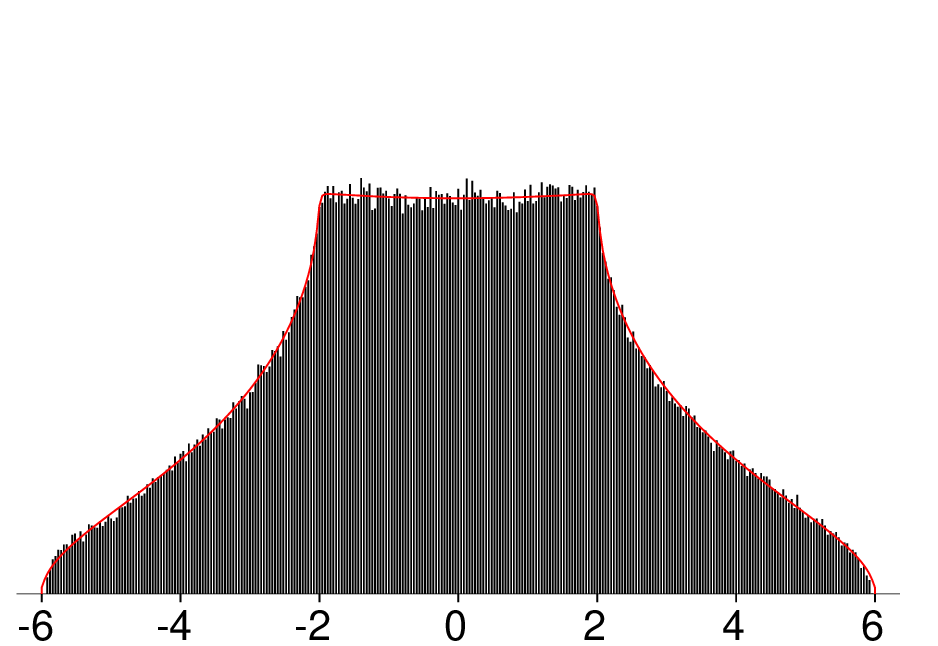}
\end{center}\vskip-5mm
\caption{Theoretical and experimental trace distributions for the
$K3$~surface
$X_5$
over a non-normal cubic number field
$k$
of geometric Picard
rank~$16$
having CM by
$\smash{E=k(i)}$
}\vskip-2mm
\end{figure}
\end{ex}

\subsubsection*{Real multiplication}\leavevmode\smallskip

\noindent
In the generic and RM cases, it may easily happen that the splitting field of
$\smash{G_{X,\overline{l}}/G_{X,\overline{l}}^0}$
strictly
contains~$kE$.
Cf.\ \cite[Examples 5.4, 5.8 and~5.10]{EJ22b}. In~the example below, however, equality holds, at least~conjecturally.

\begin{ex}
\label{Qw3_twist}
Let~$X_6'$
be the double cover
of~$\Pb^2_\bbQ$,
given by
$$w^2 = -1974 xyz (x^3 - 14x^2z + 11xy^2 - xz^2 + 12y^3 - 14y^2z - 12yz^2 + 14z^3)$$
and
$X_6$
the
$K3$~surface
obtained as the minimal desingularisation
of~$X_6'$.

\begin{abc}
\item
Then the geometric Picard rank
of~$X_6$
is~$16$.
\item
The endomorphism field
of~$X_6$
is at most~quadratic.
\end{abc}\smallskip

\noindent
{\bf Proof.}
The~surface
$X_6$
is the quadratic twist of the
surface~$\widetilde{X}_6$,
presented in \cite[Example 5.10]{EJ22b}, by the twist
factor~\mbox{$(-1974) = -2 \!\cdot\! 3 \!\cdot\! 7 \!\cdot\! 47$}.
I.e., the two surfaces are geometrically isomorphic. In fact, they are isomorphic over
$K = \bbQ(\sqrt{-1974})$.
In~particular, the geometric Picard ranks as well as the endomorphism fields are the~same.
\eop
\end{ex}

There is strong evidence that the endomorphism field of
$X_6$
is
$\smash{E=\bbQ(\sqrt{3})}$.
Note~that
$\smash{\widetilde{X}_6 = V^{(3)}_{1,2}}$
in the notation of \cite[Conjectures~5.2]{EJ16} and
$\smash{\widetilde{X}_6 = X_6}$
in the notation of~\cite{EJ22b}. Thus,~conjecturally,
$\smash{G_{X_6,\overline{l}}^0 \cong [\SO_3(\overline\bbQ_l)]^2}$,
for any
prime~$l$.

Concerning~the component group, we claim that
$\smash{G_{X_6,\overline{l}}/G_{X_6,\overline{l}}^0}$
is of
order~$2$.
For~this, recall the~inclusion
$$G_{X_6,\overline{l}}/G_{X_6,\overline{l}}^0 \hookrightarrow \N_{\O_6(\overline\bbQ_l)}([\SO_3(\overline\bbQ_l)]^2)/[\SO_3(\overline\bbQ_l)]^2 \cong (\bbZ/2\bbZ)^2 \!\rtimes\! S_2$$
into the dihedral group of
order~$8$.
Moreover,~one clearly has
$$\#X_6'(\bbF_{\!p}) = 
\left\{
\begin{array}{ll}
\#\widetilde{X}'_6(\bbF_{\!p}) \,,              & {\rm ~if~} (\frac{-1974}p) = 1 \,, \\
2(p^2+p+1) - \#\widetilde{X}'_6(\bbF_{\!p}) \,, & {\rm ~if~} (\frac{-1974}p) = -1
\end{array}
\right.
$$
and, consequently,
\begin{equation}
\label{rep_twist}
\textstyle \varrho_{X_6,\overline{l}} = \varrho_{\widetilde{X}_6,\overline{l}} \!\cdot\! (\frac{-1974}\cdot) \,.
\end{equation}
Thus,~the experimental observations described in \cite[Section~5]{EJ22b} carry over as~follows.

One~has
$\smash{\Tr(\varrho_{X_6,\overline{l}}(\Frob_p)) = 0}$
for all primes
$p = \pm5 \pmod {12}$,
at least as long
as~\mbox{$p < 5 \!\cdot\! 10^8$}.
These~are exactly the primes such that the jump character
$\tau_{X_6}$
evaluates
to~$(-1)$
at the corresponding Frobenii, hence the component group is bound to elements of the types
$\smash{\genfrac(){0pt}{}{+ \,\,0}{\,0 \,\,+}}$,
$\smash{\genfrac(){0pt}{}{- \,\,0}{\,0 \,\,-}}$,
$\smash{\genfrac(){0pt}{}{\,0 \,\,+}{+ \,\,0}}$,
and
$\smash{\genfrac(){0pt}{}{\,0 \,\,-}{- \,\,0}}$.
Furthermore, formula~(\ref{rep_twist}) shows that, in view of the corresponding \mbox{observation}
for~$\smash{\widetilde{X}_6}$
\cite[comments below Example~5.10]{EJ22b},
$$\textstyle \varrho_{X_6,\overline{l}}(\Frob_p) \in [\O_3^-(\overline\bbQ_l)]^2 \;\,\Longleftrightarrow\;\; p \equiv \pm1 \pmodulo{12} {\rm ~and~} (\frac3p) = -1 \,,$$
which is~contradictory. Hence,~the component group
$\smash{G_{X_6,\overline{l}}/G_{X_6,\overline{l}}^0}$
is bound to
order~$2$.

The~experimental distribution of the traces for all primes
$p < 5\cdot10^8$
is shown in the histogram below, plotted against the theoretical distribution, according to the Sato--Tate conjecture, the density of which is given~by
\[
\smash{\textstyle \frac1{8\pi^2} \big((2-t) \K(1-\frac{(t-2)^2}{16}) + 4 \E(1-\frac{(t-2)^2}{16}) \big) \,,}
\]
for
$\K$
and~$\E$
the complete elliptic integrals of the first and second~kinds. Cf.\ \cite[Section~3]{EJ22b}. The~spike of
mass~$1/2$
above~$0$
is not~shown.\vskip-10mm

\begin{figure}[H]
\label{RM_plot}
\begin{center}
\includegraphics[scale=0.4]{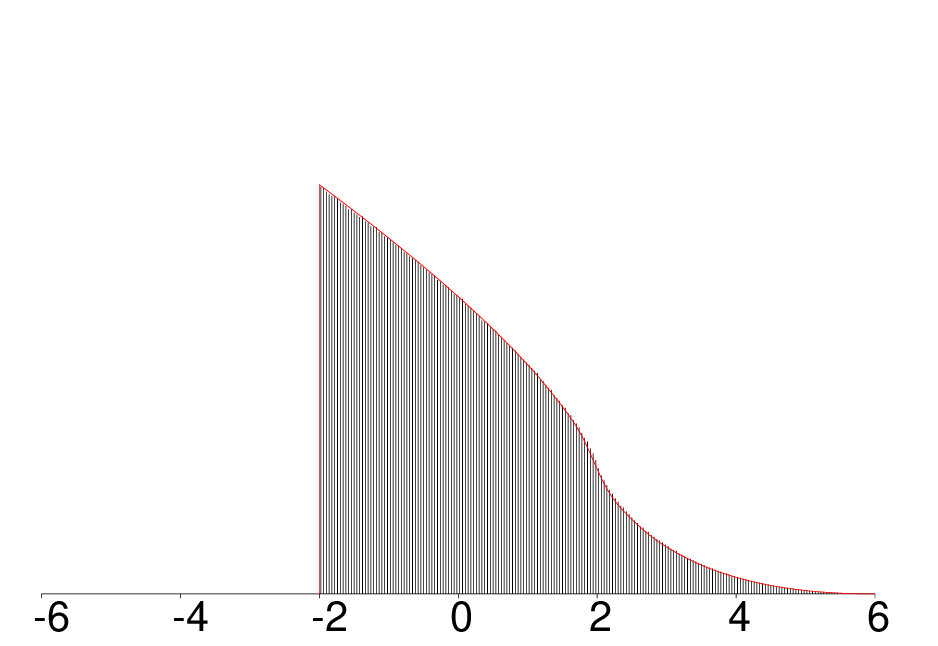}
\end{center}\vskip-5mm
\caption{Theoretical and experimental trace distributions for the
$K3$~surface
$X_6$
over
$k=\bbQ$
of geometric Picard
rank~$16$
having RM by
$\smash{E=\bbQ(\sqrt{3})}$
}
\end{figure}\vskip-3mm

\noindent
\begin{rem}
Put~$X_{6'} := (X_6)_k$,
for~$\smash{k = \bbQ(\sqrt{3})}$.
Then~the geometric Picard rank and the endomorphism
field~$E$
are the same as
for~$X_6$.
However,~if
$\smash{E = \bbQ(\sqrt{3})}$
is indeed true then
$\smash{G_{X_{6'},\overline{l}} = [\SO_3(\overline\bbQ_l)]^2}$
is~connected. I.e.,~in the histogram above, the spike disappears. As~shown in Theorem~\ref{main}, this happens if and only
if~$k \supseteq E$.
\end{rem}

\setlength\parindent{0mm}
\end{document}